\documentclass[reqno,draft]{amsart}
\usepackage{amsmath}
\allowdisplaybreaks[3]
\numberwithin{equation}{section}
\numberwithin{equation}{section}

\def\proof{\indent{\em Proof.\quad}}
\def\endproof{\hfill\hbox{$\sqcup$}\llap{\hbox{$\sqcap$}}\medskip}
\newtheorem{thm}{{\indent\bf Theorem}}[section]
\newtheorem{prop}{{\indent\bf Proposition}}[section]
\newtheorem{dfn}{{\indent\bf Definition}}[section]
\newtheorem{rmk}{{\indent\bf Remark}}[section]
\newtheorem{lem}{{\indent\bf Lemma}}[section]
\newtheorem{cor}{{\indent\bf Corollary}}[section]
\newtheorem{expl}{{\indent\bf Example}}[section]

\newcommand{\mb}{\mbox}
\newcommand{\hs}{\hspace}
\newcommand{\dps}{\displaystyle}

\newcommand{\vs}{\vspace}
\newcommand{\ol}{\overline}

\newcommand{\strl}{\stackrel}

\newcommand{\td}{\tilde}

\newcommand{\fr}{\frac}
\newcommand{\ed}{{\rm End}}
\newcommand{\edd}{\end{document}}
\newcommand{\be}{\begin{equation}}
\newcommand{\ee}{\end{equation}}

\newcommand{\lagl}{\langle}
\newcommand{\ragl}{\rangle}
\newcommand{\lmx}{\left(\begin{matrix}}
\newcommand{\rmx}{\end{matrix}\right)}
\newcommand{\ldt}{\left|\begin{matrix}}
\newcommand{\rdt}{\end{matrix}\right|}

\newcommand{\tr}{{\rm tr\,}}
\newcommand{\vfi}{\varphi}

\newcommand{\bbr}{{\mathbb R}}

\newcommand{\bbc}{{\mathbb C}}
\newcommand{\bbq}{{\mathbb Q}}
\newcommand{\ba}{\begin{array}}
\newcommand{\ea}{\end{array}}
\newcommand{\nnm}{\nonumber}
\newcommand{\beal}{\begin{align}}
\newcommand{\eal}{\end{align}}
\newcommand{\bea}{\begin{eqnarray}}
\newcommand{\eea}{\end{eqnarray}}

\newcommand{\spn}{{\rm Span\,}}

\newcommand{\pp}[2]{\fr{\partial #1}{\partial #2}}
\newcommand{\ppp}[3]{\fr{\partial^2 #1}{\partial #2\partial #3}}
\newcommand{\dd}[2]{\fr{d #1}{d #2}}
\begin{document}

\title[The symmetric equiaffine hyperspheres and
the symmetric Lagrangian submanifolds]{On the symmetric equiaffine hyperspheres and\\
the minimal symmetric Lagrangian submanifolds}%
\author{Xingxiao Li}%

\date{}

\begin{abstract}
In this paper, a correspondence via duality is established between the set of locally strongly convex symmetric equiaffine hyperspheres and the set of minimal symmetric Lagrangian submanifolds in a certain complex space form. By using this correspondence theorem, we are able to provide an alternative proof of the classification theorem for the locally strongly convex equiaffine hypersurfaces with parallel Fubini-Pick forms, which has been established recently by Z.J. Hu etc in a totally different way.
\end{abstract}
\maketitle

\tableofcontents

\section{Introduction}

In differential geometry there are many different research branches of interest and among them are two well known important ones: the differential geometry of Lagrangian submanifolds in complex space forms and the affine differential geometry of nondegenerate hypersurfaces.
As we know, finding the links or relations between different research branches of mathematics, in particular those of differential geometry, is of great interest and importance for us. In this article we are to start the consideration of the relation between the differential geometry of Lagrangian submanifolds in complex space forms and the equiaffine differential geometry of nondegenerate hypersurfaces. In fact, we presently first consider the Riemannian case. The more general case for the pseudo-Riemannian case will be considered in some forthcoming papers.

Lagrangian submanifolds of complex space forms are very special and interesting class of objects in the theory of submanifolds and have been studied extensively. Recent years various kinds of research achievement on this subject have been obtained. In particular, the study of minimal Lagrangian submanifolds seems more interesting and more attracting. For example, the classification of minimal Lagrangian submanifolds of constant curvatures in the complex projective space $\bbc P^n$ (\cite{li-zhao94}), that of parallel Lagrangian submanifolds in complex space forms (\cite{nai81}, \cite{nai83a}, \cite{nai83b}, \cite{nai-tak82}), etc.

On the other hand, affine hypersperes are very special in the equiaffine differential geometry of hypersurfaces. In particular, if an affine hypersurface is of parallel Fubini-Pick form, then it must be an affine hypersphere (\cite{bok-nom-sim90}). The study of affine hypersurfaces is also fruitful in recent ten years. As we know, affine hyperspheres seems simple in definition but they do form a very large class of hypersurfaces. Thus to find all the affine hyperspheres has been a great challenge and still remains a very hard job to be done. Although this, the study of affine hyperspheres have been made a lot of great achievment by many authors. For example, the proof of the Calabi's conjecture (\cite{amli90}, \cite{amli92}), the classification of hyperspheres of constant affine curvatures (\cite{vra-li-sim91},  \cite{wang93} and \cite{kri-vra99}), the generalizations of Calabi's composition of affine hyperbolic hyperspheres (with multiple factors, \cite{lix93}; for more general cases, \cite{dil-vra94}), the characterization of the Calabi's composition of hyperbolic hyperspheres (\cite{hu-li-vra08}), and the classification of locally strongly convex hypersurfaces with parallel Fubini-Pick forms (for some special cases, \cite{dil-vra-yap94}, \cite{hu-li-sim-vra09}, complete for general case, \cite{hu-li-vra11}). As for the general nondegenerate case, there also have been some interesting partial classification results, see for example the series of published papers by Z.J. Hu etc: \cite{hu-li11}, \cite{hu-li-li-vra11a} and \cite{hu-li-li-vra11b}. In this direction, a very recent development is the preprint article \cite{Hil12} in which the author aimed at a complete classification of nondegenerate centroaffine hypersurfaces with parallel Fubini-Pick form.

We particularly remark that, F. Dillen, H.Z. Li and X.F. Wang introduced and studied the Calabi's composition of parallel Lagrangian submanifolds in the complex projective space $\bbc P^n$ (\cite{dil-li-wang12}); Recently, by using the idea and techniques developed in \cite{hu-li-sim-vra09} and \cite{li-wang11}, H.Z. Li and X.F. Wang also gave a totally different proof (in fact in a geometric manner) of the complete classification of all parallel Lagrangian submanifolds in $\bbc P^n$ (\cite{dil-li-wang12}).

In this paper, we find a close link between the locally strongly convex symmetric equiaffine hyperspheres and the minimal and symmetric Lagrangian submanifolds in the complex space forms. By observing the apparent similarity between the Gaussian equations of the equiaffine hyperspheres and the minimal Lagrangian submanifolds in the complex space forms, we use the duality of Riemannian symmetric spaces in order to establish a direct correspondence between the set of affine equivalent class of locally strongly convex symmetric affine hyperspheres and the set of holomorphic isometric class of minimal symmetric Lagrangian submanifolds in a certain complex space form. See Theorem \ref{corr thm}. By making use of this correspondence theorem, we are able to provide an alternative proof of the classification theorem (Theorem \ref{cla thm}) for the locally strongly convex equiaffine hypersurfaces with parallel Fubini-Pick forms, which has been established recently by Z.J. Hu etc in a totally different way (see \cite{hu-li-vra11} for the detail).

{\sc Acknowledgement} The author is grateful to Professor A-M Li for his encouragement and important suggestions during the preparation of this article. He also thanks Professor Z.J. Hu for providing him valuable related references some of which are listed in the end of this paper. The main results of this paper has been announced at the conference of differential geometry, June, 2012, Bedlewo, Poland. The author would like to express his hearty thanks to the conference organizers, in particular, Professor Barbara Opozda and Professor Udo Simon for their kind invitation and hospitality.

\section{Preliminaries}

\subsection{The equiaffine differential geometry of hypersurfaces}

Let $x:M^n\to\bbr^{n+1}$ be nondegenerate hypersurface. Then there are several basic equiaffine invariants of $x$ among which are: the affine metric (Berwald-Blaschke metric) $g$, the affine normal $\xi:=\fr1n\Delta_gx$, the Fubini-Pick $3$-form (the so called cubic form) $A\in\bigodot^3T^*M^n$ and the affine second fundamental $2$-form $B\in\bigodot^2T^*M^n$. By using the index lifting by the metric $g$, we can identify $A$ and $B$ with the linear maps $A:TM\to \ed(TM)$ or $A:TM\bigodot TM\to TM$ and $B:TM\to TM$, respectively, by
\be\label{ab}
g(A(X)Y,Z)=A(X,Y,Z) \mb{\ or\ }g(A(X,Y),Z)=A(X,Y,Z),\quad
g(B(X),Y)=B(X,Y),
\ee
for all $X,Y,Z\in TM$. Sometimes we call the corresponding $B\in \ed(TM)$ the affine shape operator of $x$. In this sense, the affine Gauss equation can be written as follows:
\be\label{gaus}
R(X,Y)Z=\fr12(g(Y,Z)B(X)+B(Y,Z)X-g(X,Z)B(Y)-B(X,Z)Y)-[A(X),A(Y)](Z),
\ee
where, for any linear transformations $T,S\in \ed(TM)$,
\be\label{comm}
[T,S]=T\circ S-S\circ T.
\ee
Each of the eigenvalues $B_1,\cdots,B_n$ of the linear map $B:TM\to TM$ is called the affine principal curvature of $x$. Define
\be\label{afme}
L_1:=\fr1n\tr B=\fr1n\sum_iB_i.
\ee
Then $L_1$ is referred to as the affine mean curvature of $x$. A hypersurface $x$ is called an (elliptic, parabolic, or hyperbolic) affine hypersphere, if all of its affine principal curvatures are equal to one (positive, 0, or negative) constant. In this case we have
\be\label{afsp}
B(X)=L_1X,\quad\mb{for all\ }X\in TM.
\ee
It follows that the affine Gauss equation \eqref{gaus} of an affine hypersphere assumes the following form:
\be\label{gaus_af sph}
R(X,Y)Z=L_1(g(Y,Z)X-g(X,Z)Y)-[A(X),A(Y)](Z),
\ee

Furthermore, all the affine lines of an elliptic affine hypersphere or a hyperbolic affine hypersphere $x:M^n\to\bbr^{n+1}$ pass through a fix point $o$ which is refer to as the affine center of $x$; Both the elliptic affine hyperspheres and the hyperbolic affine hyperspheres are called proper affine hyperspheres, while the parabolic affine hyperspheres are called improper affine hyperspheres.

{\prop\label{affine spheres} $($\cite{li-sim-zhao93}$)$ A nondegenerate hypersurface $x:M^n\to \bbr^{n+1}$ is a proper affine hypersphere with affine mean curvature $L_1$ and with the origin $o$ as its affine center if and only if the affine line is parallel to the position vector $x$. In this case, the affine normal $\xi$ is given by $\xi=-L_1x$.}

For each vector field $\eta$ transversal to the tangent space of $x$, we have the following direct decomposition
$$
x^*T\bbr^{n+1}=x_*(TM)\oplus \bbr\cdot\eta.
$$
This decomposition and the canonical differentiation $\bar D^0$ on $\bbr^{n+1}$ define a bilinear form $h\in\bigodot^2T^*M^n$ and a connection $D^\eta$ on $TM$ as follows:
\be\label{dfn h}
\bar D^0_XY=x_*(D^\eta_XY)+h(X,Y)\eta,\quad\forall X,Y\in TM.
\ee
\eqref{dfn h} can be referred as to the {\em affine Gauss formula} of the hypersurface $x$.
In particular, in case that $\eta$ is parallel to the affine normal $\xi$, the induced connection $\nabla:=D^\eta$ is independent of the choice of $\eta$ and is referred to as the affine connection of $x$.

In what follows we make the following convention for the range of indices:
$$1\leq i,j,k,l\leq n.$$

Let $\{e_i,e_{n+1}\}$ be a local unimodular frame field along $x$ with $e_{n+1}$ parallel to the affine normal $\xi$, and $\{\omega^i,\omega^{n+1}\}$ be its dual coframe. Then the above invariants can be respectively expressed locally as
\be\label{gab}
g=\sum g_{ij}\omega^i\omega^j,\quad A=\sum A_{ijk}\omega^i\omega^j\omega^k,\quad B=\sum B_{ij}\omega^i\omega^j.
\ee

Now we directly write down the basic formulas of the equiaffine geometry for the affine hypersurface $x$. For the details, please see \cite{li-sim-zhao93} and \cite{nom-sas94}.
\begin{align}
&\sum_{i,j} g^{ij}A_{ijk}=0,\label{basic1}\\
&R_{ijkl}=\sum_m(A^m_{ik}A_{jlm}-A^m_{il}A_{jkm})+\fr12(g_{il}B_{jk}+g_{jk}B_{il} -g_{ik}B_{jl}-g_{jl}B_{ik}),\label{basic2}\\
&R_{ij}=\sum_{k,l}A^k_{il}A^l_{jk}+\fr n2L_1g_{ij}+\fr{n-2}2 B_{ij},\label{basic2-1}\\
&A_{ijk,l}-A_{ijl,k}=\fr12(g_{ik}B_{jl}+g_{jl}B_{ik} -g_{il}B_{jk}-g_{jk}B_{il}),\label{basic3}\\
&\sum_{l}A^l_{ij,l}=\fr n2(L_1g_{ij}-B_{ij}),\label{basic3-1}
\end{align}
where $R_{ijkl}$ are the components of the Riemannian curvature tensor of the Berwald-Blaschke metric $g$,
while $A_{ijk,l}$ and $A_{ijk,lm}$ are the covariant derivatives of $A_{ijk}$ with respect to Levi-Civita connection of $g$.

Define the normalized scalar curvature $\chi$ and the Pick invariant $J$ by
$$
\chi=\fr1{n(n-1)}\sum g^{il}g^{jk}R_{ijkl},\quad J=\fr1{n(n-1)}\sum A_{ijk}A_{pqr}g^{ip}g^{jq}g^{kr}.$$
Then the affine Gauss equation can be written in terms of the metric and the Fubini-Pick form as follows
\begin{align}
R_{ijkl}=&(A_{ijk,l}-A_{ijl,k})+(\chi-J)(g_{il}g_{jk}-g_{ik}g_{jl})\nnm\\ &\ +\fr2n\sum(g_{ik}A_{jlm,m}-g_{il}A_{jkm,m}) +\sum_m(A^m_{ik}A_{jlm}-A^m_{il}A_{jkm}).
\label{basic2'}\end{align}
Write $h=\sum h_{ij}\omega^i\omega^j$ and $H=\det(h_{ij})$. Then
\be\label{dfn g}
g_{ij}=H^{-\fr1{n+2}}h_{ij},\quad \xi=H^{\fr1{n+2}}e_{n+1}.
\ee
Define
\be\label{hijk0}
\sum_kh_{ijk}\omega^k=dh_{ij}+h_{ij}\omega^{n+1}_{n+1}-\sum h_{kj}\omega^k_i-\sum h_{ik}\omega^k_j.
\ee
Then the Fubini-Pick form $A$ can be determined by the following formula:
\be\label{hijktoaijk}
A_{ijk}=-\fr12H^{-\fr1{n+2}}h_{ijk}.
\ee

{\dfn A nondegenerate hypersurface $x:M^n\to \bbr^{n+1}$ is called affine symmetric (resp. locally affine symmetric) if

(1) the pseudo-Riemannian manifold $(M^n,g)$ is symmetric (resp. locally symmetric) and therefore $(M^n,g)$ can be written (resp. locally written) as $G/K$ for some connected Lie group $G$ of isometries with $K$ one of its closed subgroups;

(2) the Fubini-Pick form $A$ is invariant under the action of $G$.}

{\prop\label{sym para1}
A nondegenerate hypersurface $x:M^n\to\bbr^{n+1}$ is of parallel Fubini-Pick form $A$ if and only if $x$ is locally affine symmetric.}

\proof
First we suppose that the Fubini-Pick form $A$ of $x$ is parallel. Then by \cite{bok-nom-sim90}, $x$ must be an affine hypersphere. It then follows from \eqref{gaus_af sph} that the Berwald-Blaschke metric $g$ must be locally symmetric. Thus locally we can write $M^n=G/K$ and the canonical decomposition of the corresponding orthogonal symmetric pair $(\mathfrak{g},\mathfrak{k})$ is written as
${\mathfrak g}=\mathfrak{k}+\mathfrak{m}$ where the vector space $\mathfrak m$ is identified with $T_oM$. Here $o\in M^n$ is the base point given by $o=eK$ with $e$ the identity of $G$.
Note that, for all $X,Y_i\in {\mathfrak m}=T_oM$, $i=1,2,3$, the vector field $Y_i(t):=L_{\exp (tX)*}(Y_i)$ is the parallel translation of $Y_i$ along the geodesic $\gamma(t):=_{\exp (tX)}\!\!K$ (see, for example, \cite{hel01}).
Consequently we have
\begin{align}
&\dd{}{t}((L_{\exp (tX)}^*A)(Y_1,Y_2,Y_3))\nnm\\
=&\dd{}{t}(A_{\exp (tX)K}(L_{\exp (tX)*}(Y_1),
L_{\exp (tX)*}(Y_2),L_{\exp (tX)*}(Y_3)))\nnm\\
=&(\hat\nabla_{\gamma'(t)}A)(Y_1(t),Y_2(t),Y_3(t))=0,
\end{align}
where $\hat{\nabla}$ is the Levi-Civita connection of the metric $g$.
It follows that
\be\label{2.19-0}
A_{\exp (tX)K}(L_{\exp (tX)*}(Y_1),
L_{\exp (tX)*}(Y_2),L_{\exp (tX)*}(Y_3))
\ee
is constant with respect to the parameter $t$ and thus $A$ is $G$-invariant.

Conversely, we suppose that $M^n=G/K$ locally for some symmetric pair $(G,K)$ and that $A$ is $G$-invariant. Then for any $X,Y_i\in {\mathfrak m}=T_oM$, $i=1,2,3$, the function
\eqref{2.19-0}
is again a constant along the geodesic $\gamma(t)$.

Therefore,
$$
(\hat\nabla_XA)(Y_1,Y_2,Y_3)=\left.\dd{}{t}\right|_{t=0} A_{\gamma(t)}(Y_1(t),Y_2(t),Y_3(t))=0,
$$
where we have once again used the fact that each $Y_i(t)$ is parallel along the geodesic $\gamma(t)$.
\endproof

The following existence and uniqueness theorems are well known:

{\thm\label{affine existence} $($\cite{li-sim-zhao93}$)$
$($The existence$)$ Let $(M^n,g)$ be a simply connected Riemannian manifold
of dimension $n$, and $A$ be a symmetric $3$-form on $M^n$ satisfying the
affine Gauss equation \eqref{basic2} and the apolarity condition \eqref{basic1}. Then there exists a locally strongly convex immersion $x:M^n\to \bbr^{n+1}$ such that $g$ and $A$ are the Berwald-Blaschke metric and the Fubini-Pick form for $x$, respectively.}

{\thm\label{affine uniqueness} $($\cite{li-sim-zhao93}$)$ $($The uniqueness$)$ Let $x:M^n\to \bbr^{n+1}$,
$\bar x:\bar M^n\to \bbr^{n+1}$ be two locally strongly convex hypersurfaces of dimension $n$ with respectively the Berwald-Blaschke metrics $g$, $\bar g$ and the Fubini-Pick forms $A$, $\bar A$, and $\vfi:(M^n,g)\to (\bar M^n,\bar g)$ be an isometry between Riemannian manifolds. Then $\vfi^*\bar A=A$ if and only if there exists a unimodular affine transformation $\Phi:\bbr^{n+1}\to \bbr^{n+1}$ such that $\bar x\circ\vfi=\Phi\circ x$, or equivalently, $\bar x=\Phi\circ x\circ\vfi^{-1}$.}

\rmk\rm The necessity part of Theorem \ref{affine uniqueness} is proved in \cite{li-sim-zhao93}. Here we give a proof for the sufficient part as follows:

Choose an orthonormal frame field $\{e_i;\ 1\leq i\leq n\}$ on $M^n$ with its dual coframe $\{\omega^i;\ 1\leq i\leq n\}$. Let $\xi,\bar \xi$ are respectively the affine normal of $x$ and $\bar x$. Then $\{e_1,\cdots,e_n,\xi\}$ is unimodular.
Define $\bar e_i=\vfi_*(e_i)$, $\bar\omega^i=(\vfi^{-1})^*\omega^i$, $1\leq i\leq n$. Then $\{\bar\omega^i;\ 1\leq i\leq n\}$ is the dual coframe of $\{\bar e_i;\ 1\leq i\leq n\}$. Since $\vfi$ is an isometry, $\{\bar e_1,\cdots,\bar e_n,\bar \xi\}$ is also unimodular.

Under the condition that $\bar x=\Phi\circ x\circ\vfi^{-1}$, we claim that $\bar \xi=(\Phi_*(\xi))\circ {\vfi^{-1}}$. In fact
\be\label{bar ejei}
\bar e_j(\bar e_i\bar x)=\vfi_*(e_j)(\vfi_*(e_i)(\Phi\circ x\circ\vfi^{-1})) =\vfi_*(e_j)
((e_i(\Phi\circ x))\circ {\vfi^{-1}})
=(e_j(e_i(\Phi\circ x)))\circ {\vfi^{-1}}.
\ee
Denote respectively by $\nabla,\hat\nabla,\Delta$ and $\bar\nabla,\hat{\bar\nabla},\bar\Delta$ the affine connections of $x,\bar x$, the Riemannian connections and the Laplacians of $g,\bar g$. Then we find
\begin{align}
\bar \xi=&\fr1n\bar\Delta\bar x=\fr1n\left(\sum_i\left(\bar e_i(\bar e_i\bar x)-(\hat{\bar\nabla}_{\bar e_i}\bar e_i)(x)\right)\right)\nnm\\
=&\fr1n\left(\sum_i\left((e_i(e_i(\Phi\circ x)))\circ {\vfi^{-1}}-\vfi_*(\hat\nabla_{e_i}e_i)(\bar x)\right)\right)\nnm\\
=&\fr1n\left(\sum_i\left((e_i(\Phi_*(e_ix)))\circ {\vfi^{-1}} -(\hat\nabla_{e_i}e_i)(\Phi\circ x)\circ{\vfi^{-1}}\right)\right)\nnm\\
=&\fr1n\left(\sum_i\left((\Phi_*(e_i(e_ix)))\circ {\vfi^{-1}} -(\Phi_*(\hat\nabla_{e_i}e_i)(x))\circ {\vfi^{-1}}\right)\right)\nnm\\
=&\fr1n\Phi_*\left(\sum_i\left(e_i(e_ix) -(\hat\nabla_{e_i}e_i)(x)\right)\right)\circ {\vfi^{-1}} \nnm\\
=&\fr1n(\Phi_*(\Delta x))\circ {\vfi^{-1}}
=(\Phi_*(\xi))\circ {\vfi^{-1}}.\nnm
\end{align}

On the other hand, by \eqref{bar ejei} and the affine Gauss formula \eqref{dfn h} of $x$
\begin{align}
\bar e_j\bar e_i \bar x=&(e_je_i(\Phi\circ x))\circ\vfi^{-1} =(e_j\Phi_*(e_i(x))\circ\vfi^{-1} =(\Phi_*(e_je_i(x))\circ\vfi^{-1}\nnm\\
=&\left(\Phi_*\left(x_*(\nabla_{e_j}e_i)+\delta_{ij}\xi\right)\right)\circ\vfi^{-1}\nnm\\
=&\left(\Phi_*\left(x_*(\nabla_{e_j}e_i)\right)\right)\circ\vfi^{-1} +\delta_{ij}\left(\Phi_*(\xi)\right)\circ\vfi^{-1}.\nnm
\end{align}
But, by the affine Gauss formula of $\bar x$,
$$
\bar e_j\bar e_i \bar x=\bar x_*(\bar\nabla_{\bar e_j}\bar e_i)+\delta_{ij}\bar \xi =\Phi_*\left(x_*(\vfi^{-1}_*(\bar\nabla_{\bar e_j}\bar e_i))\right)\circ\vfi^{-1} +\delta_{ij}(\Phi_*(\xi))\circ\vfi^{-1}.
$$
It follows that
$$
\Phi_*\left(x_*(\nabla_{e_j}e_i)\right)=\Phi_*\left(x_*(\vfi^{-1}_*(\bar\nabla_{\bar e_j}\bar e_i))\right).
$$
Therefore
\be
\vfi^{-1}_*(\bar\nabla_{\bar e_j}\bar e_i)=\nabla_{e_j}e_i,
\text{\ or equivalently,\ }
\vfi_*(\nabla_{e_j}e_i)=\bar\nabla_{\bar e_j}\bar e_i,
\ee
from which we find that
\begin{align}
\bar A(\bar e_i,\bar e_j,\bar e_k)=&\bar g(\bar A(\bar e_i,\bar e_j),\bar e_k)
=\bar g(\bar\nabla_{\bar e_j}\bar e_i-\hat{\bar\nabla}_{\bar e_j}\bar e_i,\bar e_k)\nnm\\
=&g(\vfi^{-1}_*(\bar\nabla_{\bar e_j}\bar e_i)-\vfi^{-1}_*(\hat{\bar\nabla}_{\bar e_j}\bar e_i),e_k) =g((\nabla_{e_j}e_i-\hat{\nabla}_{e_j}e_i),e_k) = A(e_i,e_j,e_k).
\end{align}
Consequently
\begin{align}
\bar A=&\sum \bar A(\bar e_i,\bar e_j,\bar e_k)\bar\omega^i\bar\omega^j\bar\omega^k =\sum A(e_i,e_j,e_k)(\vfi^{-1})^*\omega^i(\vfi^{-1})^*\omega^j(\vfi^{-1})^*\omega^k\nnm\\ =&(\vfi^{-1})^*\left(\sum A(e_i,e_j,e_k)\omega^i\omega^j\omega^k\right)=(\vfi^{-1})^*A,
\end{align}
or equivalently, $\vfi^*\bar A=A$. We are done.

\rmk\label{rmk2.2}\rm It is not hard to see that Theorems \ref{affine existence} and \ref{affine uniqueness} still hold for the general nondegenerate hypersurfaces.

Motivated by Theorem \ref{affine uniqueness}, we introduce the following modified equiaffine equivalence relation between nondegenerate hypersurfaces:

{\dfn Let $x:M^n\to \bbr^{n+1}$ be a nondegenerate hypersurface with the Berwald-Blaschke metric $g$. A hypersurface $\bar x:M^n\to \bbr^{n+1}$ is called affine equivalent to $x$ if there exists a unimodular transformation $\Phi:\bbr^{n+1}\to \bbr^{n+1}$ and an isometry $\vfi$ of $(M^n,g)$ such that $\bar x=\Phi\circ x\circ\vfi^{-1}$}.

For a pseudo-Riemannian manifold $(M^n,g)$, denote by $I(M^n)$ its isometry group. Given a fixed point $o\in M^n$, define a subgroup ${\mathcal F}_o(M^n)$ of isometries
$$
{\mathcal F}_o(M^n)=\{\phi\in I(M^n);\ \phi(o)=o\}.
$$
Then ${\mathcal F}_o(M^n)$ acts on $\bigodot^3(T^*_oM)$ as follows: For all $\phi\in {\mathcal F}_o(M^n)$, $T\in \bigodot^3(T^*_oM)$, and $X,Y,Z\in T_{o}M^n$,
\be\label{actoffo1}
(\phi\cdot T)(X,Y,Z):=((\phi^{-1})_o^*T)(X,Y,Z) =T((\phi^{-1})_{*o}(X),(\phi^{-1})_{*o}(Y),(\phi^{-1})_{*o}(Z));
\ee
Furthermore, if we take $T$ as a symmetric $T_oM$-valued $2$-form, then we have
\be\label{actoffo2}
(\phi\cdot T)(X,Y) :=\phi_{o*}(T((\phi^{-1})_{*o}(X),(\phi^{-1})_{*o}(Y))).
\ee

{\prop Let $x:M^n\equiv G/K\to \bbr^{n+1}$ be a nondegenerate symmetric hypesurface with the Berwald-Blaschke metric $g$. Fix one point $o\in M^n$ as the base point and denote by $A_o$ the Fubini-Pick form of $x$ at $o$. Then, a symmetric hypersurface $\bar x:M^n\to \bbr^{n+1}$ with the Fubini-Pick form $\bar A_o$ at $o$ is affine equivalent to $x$ if and only if there exists an element $\phi\in {\mathcal F}_o(M^n)$ such that $\bar A_0=\phi\cdot A_0$.}

\proof Denote by $A,\bar A$, respectively, the Fubini-Pick forms of $x$ and $\bar x$, and suppose that $\bar A_0=\phi\cdot A_0$ for some $\phi\in {\mathcal F}_o(M^n)$. Define $\td x=x\circ\phi^{-1}$ and let $\td A$ be the Fubini-Pick form of $\td x$. Then by Theorem \ref{affine uniqueness} and Remark \ref{rmk2.2}, we have $\phi^*\td A=A$, which restricting to the point $o$ gives $\td A_o=(\phi^{-1})^*A_0$, or, $\td A_0=\phi\cdot A_0=\bar A_0$. Since $\phi$ is an isometry, Proposition \ref{sym para1} shows that $\td A$ is
also $G$-invariant. This shows that $\bar A\equiv \td A$. Thus the uniqueness theorem (Theorem \ref{affine uniqueness}) assures that there exists one unimodular transformation $\Phi$ on $\bbr^{n+1}$ such that $\bar x=\Phi\circ\td x$, or $\bar x=\Phi\circ x\circ\phi^{-1}$.

Conversely suppose that $\bar x:M^n\to \bbr^{n+1}$ is affine equivalent to $x:M^n\to \bbr^{n+1}$. Then $\bar x=\Phi\circ x\circ\vfi^{-1}$ for some unimodular transformation $\Phi:\bbr^{n+1}\to \bbr^{n+1}$ and some isometry $\vfi$ of $(M^n,g)$. Since $(G/K,g)$ is symmetric and $A$ is $G$-invariant, we can choose an isometry $h\in G\subset I(M^n)$, such that

(1) $h(o)=\vfi^{-1}(o)$;

(2) $h^*(A_{h(o)})=A_o$.

Put $\phi=\vfi\circ h$. Then $\phi\in {\mathcal F}_o(M^n)$ and
\begin{align}
((\vfi^{-1})^*A_{\vfi^{-1}(o)})(X,Y,Z)=&A_{\vfi^{-1}(o)}(\vfi^{-1}_*X,\vfi^{-1}_*Y,\vfi^{-1}_*Z) =A_{h\phi^{-1}(o)}(h_*\phi^{-1}_*X,h_*\phi^{-1}_*Y,h_*\phi^{-1}_*Z)\nnm\\
=&(h^*A_{h(o)})(\phi^{-1}_*X,\phi^{-1}_*Y,\phi^{-1}_*Z) =A_o(\phi^{-1}_*X,\phi^{-1}_*Y,\phi^{-1}_*Z)\nnm\\
=&((\phi^{-1})^*A_o)(X,Y,Z)=(\phi\cdot A_o)(X,Y,Z).\nnm
\end{align}
It follows from Theorem \ref{affine uniqueness} that $\bar A_o=(\vfi^{-1})^*A_{\vfi^{-1}(o)}=\phi\cdot A_o$.
\endproof

\subsection{Symmetric Lagrangian submanifolds in the semi-Hermitian complex space forms}

Let $\bbq^n\equiv\bbq^n(4c)$ be the $n$-dimensional semi-Hermitian complex space form with constant sectional curvature $4c$ and complex structure $J$. Let $\bar g$ be the corresponding $J$-invariant metric on $\bbq^n$. An isometric immersion $\td x:M^n\to \bbq^n$ of a pseudo-Riemannian manifold $(M^n,g)$ into $\bbq^n$ is called Lagrangian if $J(\td x_*TM^n)=T^\bot M^n:=(\td x_*TM)^\bot$, or equivalently, $\td x^*\bar\omega=0$, where $\bar\omega$ is the K\"ahler form on $\bbq^n$. Furthermore, $\td x$ is called parallel if the second fundamental form $\sigma$ of $\td x$ is parallel, i.e., $\td D\sigma=0$ where the covariant differentiation $\td D$ is induced by the Levi-Civita connections $\hat\nabla$ on $M^n$ and $\bar D$ on $\bbq^n$.

Now we assume that $\td x:M^n\to \bbq^n$ is Lagrangian. By means of the metric $g$, $\bar g$ and the complex structure $J$, $\sigma$ defines a symmetric trilinear form $\td\sigma$ on $M^n$ which is also identified with a $TM$-valued symmetric bilinear form, still denoted by $\td\sigma$, such that
\be\label{td sgm}
\td\sigma(X,Y,Z)=g(\td\sigma(X,Y),Z)=\bar g(\sigma(X,Y),J\td x_*Z),\quad\forall X,Y,Z\in TM^n.
\ee
Moreover, $\td\sigma$ can also be viewed as $\ed(TM)$-valued linear map defined by
$\td\sigma(X)Y:=\td\sigma(X,Y)$ for any $X,Y\in TM$. In this sense
the Gaussian equation of $\td x$ can be written as
\be\label{gauss3}
R(X,Y)Z=c (g(Y,Z)X-g(X,Z)Y)+[\td\sigma(X),\td\sigma(Y)]Z,\quad{\rm for all\ }X,Y,Z\in TM^n.
\ee

\rmk\rm Generally, any $TM$-valued symmetric bilinear map $\Psi:TM\times TM\to TM$ is identified with a trilinear map $\Psi:TM\times TM\times TM\to \bbr$ which is not totally symmetric in general but is symmetric with respect to the first two factors. If, in addition, the corresponding trilinear map $\Psi:TM\times TM\times TM\to \bbr$ is totally symmetric, then we will call the original map $\Psi:TM\times TM\to TM$ is totally symmetric. Moreover, for each $X\in TM$, we have a linear map $\Psi(X):TM\to TM$ given by $\Psi(X)Y:=\Psi(X,Y)$, $\forall\,Y\in TM$. These identifications will be frequently used in the rest of the present paper.

{\dfn
A Lagrangian immersion $\td x: \td M^n\to \bbq^n$ with the induced metric $\td g$ is called symmetric (resp. locally symmetric) if

(1) the pseudo-Riemannian manifold $(\td M^n,\td g)$ is symmetric (resp. locally symmetric) and therefore $(\td M^n,\td g)$ can be written (resp. locally written) as $\td G/ K$ for some connected Lie group $\td G$ of isometries with $K$ one of its closed subgroups;

(2) the symmetric form $\td\sigma$ in \eqref{td sgm} induced by the second fundamental form $\sigma$ of $x$ is invariant under the action of $G$.}

In the same way as Proposition \ref{sym para1}, we can prove
{\prop\label{sym para2}
A Lagrangian isometric immersion $\td x:\td M^n\to \bbq^n$ is parallel if and only if $x$ is locally symmetric.}

From now on we assume that the Lagrangian immersion $\td x: \td M^n\to\bbq^n$ is (locally) symmetric and minimal. Then we can (locally) write $ \td M^n=\td G/ K$. Clearly in this case $K\subset {\mathcal F}_{\td o}(\td M^n)$ where $\td o:=\td e K$ with $\td e\in\td G$ being the unit element. Moreover the Lie algebra $\td{\mathfrak g}$ has a canonical decomposition $\td{\mathfrak g}={\mathfrak k}\oplus \td{\mathfrak m}$ in which ${\mathfrak k}$ is the Lie algebra of $K$, and the vector space $\td{\mathfrak m}$ is identified via the natural projection with the tangent space $T_{\td o}\td M^n$ of $\td M^n$ at the base point $\td o$.

By restriction of $\td \sigma$ to the given point $\td o$, we have an $\td{\mathfrak m}$-valued symmetric form $\td\sigma:\td{\mathfrak m}\times \td{\mathfrak m}\to \td{\mathfrak m}$ satisfying the following conditions:

(1) $\td\sigma$ is totally symmetric since $\td x$ is Lagrangian;

(2) ${\mathfrak k}\cdot\td\sigma=0$ since $\td\sigma$ is invariant under the action of $K$;

(3) $\dps\td R(X,Y)Z=c (\td g(Y,Z)X-\td g(X,Z)Y+[\td\sigma(X),\td\sigma(Y)](Z)$ for all $X,Y,Z\in \td{\mathfrak m}$, where $\td R(X,Y)$ is the curvature operator of the metric $\td g$;

(4) $\tr\,\td\sigma=0$ since $\td x$ is minimal.

Let $c$ be a given constant. For any pseudo-Riemannian symmetric space $\td M^n=\td G/K$ with $\td o=\td eK\in\td M^n$
being the base point, we denote by ${\mathcal S}_{\td M^n}(c)$ (resp. $\ol{\mathcal S}_{\td M^n}(c)$) the set of
all $\td{\mathfrak m}$-valued symmetric bilinear forms $\td\sigma$, or equivalently, the corresponding $3$-forms on $\td{\mathfrak m}$, satisfying the above
conditions (1) through (4) (resp. satisfying the above
conditions (1) through (3)). It is not hard to see that the action of ${\mathcal F}_{\td o}(\td M^n)$ given by \eqref{actoffo1} or \eqref{actoffo2} on $\bigodot^2(T^*_{\td o}\td M^n)\bigotimes (T_{*\td o}\td M^n)$ or $\bigodot^3(T^*_{\td o}\td M^n)$ keeps both ${\mathcal S}_{\td M^n}(c)$ and $\ol{\mathcal S}_{\td M^n}(c)$ invariant and thus induces an action by restriction on ${\mathcal S}_{\td M^n}(c)$ (resp. $\ol{\mathcal S}_{\td M^n}(c)$).

To simplify the statement we introduce the following equivalence relations:

{\dfn Two forms $\td\sigma_1,\td\sigma_2\in {\mathcal S}_{\td M^n}(c)$ $($resp. $\td\sigma_1,\td\sigma_2\in \ol{\mathcal S}_{\td M^n}(c))$ are called equivalent to each other if they are in the same orbit under the action of ${\mathcal F}_{\td o}(\td M^n)$.}

{\dfn Let $\td x:\td M^n\to \bbq^n$ be a Lagrangian submanifold with the induced metric $\td g$. A Lagrangian submanifold $\bar{\td x}:{\td M^n}\to \bbq^n$ is called holomorphically isometric to $\td x$ if there exist a holomorphic isometry $\Psi:\bbq^n\to \bbq^n$ and an isometry $\vfi$ of $(\td M^n,\td g)$ such that $\bar{\td x}=\Psi\circ \td x\circ\vfi^{-1}$}.

In this paper we mainly consider the case of positive definite metrics. Following \cite{nai81}, we denote by $\td\sigma_{\td x}$ the totally symmetric form defined by the second fundamental form, valued at the base point $\td o$, of a Lagrangian immersion $\td x:\td M\to \bbc P^n(4c)$. Then, using Proposition \ref{sym para2} above and Theorem 2.3, Lemmas 3.2 and 3.3 in \cite{nai81}, we can obtain the following conclusion:

{\prop\label{nai81}
Let $\td M^n=\td G/K$ be a simply connected Riemannian symmetric space with ${\mathcal S}_{\td M^n}(c)\neq\emptyset$. Then for each $\td\sigma\in{\mathcal S}_{\td M^n}(c)$, there uniquely exists one minimal symmetric Lagrangian immersion $\td x:\td M^n\to \bbc P^n(4c)$ such that $\td\sigma_{\td x}=\td\sigma$.

Furthermore, two minimal symmetric Lagrangian immersion $\td x_1,\td x_2:\td M^n\to \bbc P^n(4c)$ corresponding to some given $\td \sigma_1,\td \sigma_2\in {\mathcal S}_{\td M^n}(c)$ are holomorphic isometric if and only if $\td \sigma_1,\td \sigma_2$ are equivalent, that is, there exists some $\phi\in {\mathcal F}_{\td o}(\td M^n)$ such that $\td \sigma_2=\phi\cdot\td\sigma_1$.}

\proof The first conclusion (the existence) is that of Theorem 2.3 in \cite{nai81}; The necessary part of the second conclusion is exactly the conclusion of Lemma 3.3 in \cite{nai81}. Now suppose that $\td\sigma_1=\phi\cdot\td\sigma_2$ for some $\phi\in{\mathcal F}_{\td o}(\td M^n)$. We consider the composition $\td x_2\circ\phi^{-1}$. By Lemma 3.3 in \cite{nai81}, the symmetric form defined by the second fundamental form of the composed Lagrangian immersion $\td x_2\circ\phi^{-1}$ is $\td\sigma_{\td x_2\circ\phi^{-1}}=\phi\cdot\td\sigma_2=\td\sigma_1$. Choose an element $\Phi\in {\rm SU}(n+1)$ such that
$$
\td x_1(\td o)=\Phi(\td x_2(\td o)),\quad (\td x_1)_{*\td o}=\Phi_*\circ(\td x_2)_{*\td o}\circ(\phi^{-1})_{*\td o}.
$$
It follows that
\begin{align}
&\td x_1(\td o)=(\Phi\circ\td x_2\circ\phi^{-1})(\td o),\quad (\td x_1)_{*\td o}=(\Phi\circ\td x_2\circ\phi^{-1})_{*\td o},\nnm\\
&J\circ(\td x_1)_{*\td o}=J\circ (\Phi\circ\td x_2\circ\phi^{-1})_{*\td o} =\Phi_*\circ(J\circ(\td x_2)_{*\td o}\circ(\phi^{-1})_{*\td o}).\nnm
\end{align}
By the Lemma 3.3 in \cite{nai81},
$$
\td\sigma_{\Phi\circ x_2\circ \phi^{-1}}=\phi\cdot\td\sigma_{\td x_2}=\td\sigma_{\td x_1}.
$$
Then an application of Lemma 3.2 in \cite{nai81} shows that $\Phi\circ\td x_2\circ\phi^{-1}=\td x_1$. \endproof

\subsection{The multiple Calabi product of hyperbolic affine hyperspheres}

In 1972, E. Calabi \cite{cal72} found a composition formula by which we can construct new hyperbolic affine hyperspheres from any two given ones. The present author has generalized Calabi construction to the case of multiple factors (See \cite{lix93}). Later in 1994 F. Dillen and L. Vrancken \cite{dil-vra94} generalized Calabi original composition to any two proper affine hyperspheres and gave a detailed study of these composed affine hyperspheres. They also mentioned that their construction applies to the case of multiple factors but with no details of it. For later use, we shall first make a review of some of these facts with the emphasis on the general case of multiple factors, which seems not to have appeared in the literature other than \cite{lix93}. A detailed discussion of the formulas in this section has been given in the preprint \cite{lix11}.

\newcommand{\stx}[2]{\strl{(#1)}{#2}}
\newcommand{\spec}[1]{\prod_{#1=1}^K\fr{c_{#1}^{n_{#1}+1}H_{(#1)}^{\fr1{n_{#1}+2}}} {(n_{#1}+1)(-\!\!\stx{#1}{L}_1)}}
\newcommand{\la}{\stx{a}{L}\!\!_1{}}\newcommand{\lb}{\stx{b}{L}\!\!_1{}}
\newcommand{\lc}{\stx{c}{L}\!\!_1{}}\newcommand{\lalp}{\stx{\alpha}{L}\!\!_1{}}
\newcommand{\ha}{\!\stx{a}{h}{}\!\!} 
\newcommand{\Ha}{H_{(a)}}\newcommand{\Hb}{H_{(b)}}\newcommand{\Hc}{H_{(c)}}
\newcommand{\ga}{\!\!\stx{a}{g}{}\!\!\!}\newcommand{\Ga}{\stx{a}{G}\!\!{}}
\newcommand{\gb}{\!\!\stx{b}{g}{}\!\!\!}\newcommand{\Gb}{\stx{b}{G}\!\!{}}
\newcommand{\galp}{\!\!\stx{\alpha}{g}{}\!\!\!}
\newcommand{\xai}{x_{a,i_a}} \newcommand{\xaij}{x_{a,i_aj_a}}
\newcommand{\HH}{f_K\prod_a\fr{c_a^{(n_a+1)(f_K-1)}\Ha^{\fr{f_K+1}{n_a+2}}}
{(n_a+1)^{f_K-n_a}(-\!\!\la)^{f_K-n_a-1}}}
\newcommand{\h}{f^{-\fr1{n+2}}_K \prod_a\fr {(n_a+1)^{\fr{f_K-n_a}{f_K+1}}(-\!\!\la)^{\fr{f_K-n_a-1}{f_K+1}}}
{\left(c_a^{n_a+1}\right)^{\fr{f_K-1}{f_K+1}}\Ha^{\fr1{n_a+2}}}}
\newcommand{\oma}{\stx{a}{\omega}{}\!\!}
\newcommand{\omb}{\stx{b}{\omega}{}\!\!}
\newcommand{\tdca}{(n_a+1)\big(-\!\!\la\big)\prod_b\fr{c_b^{n_b+1}} {(n_b+1)\big(-\!\!\lb\big)}}
\newcommand{\cha}{\prod_a\fr{c_a^{n_a+1}\Ha^{\fr1{n_a+2}}} {(n_a+1)\big(-\!\!\la\big)}}
\newcommand{\chb}{\prod_b\fr{c_b^{n_b+1}\Hb^{\fr1{n_b+2}}} {(n_b+1)\big(-\!\!\lb\big)}}
\newcommand{\chc}{\prod_c\fr{c_c^{n_c+1}\Hc^{\fr1{n_c+2}}} {(n_c+1)\big(-\!\!\lc\big)}}
\newcommand{\ca}{\prod_a\fr{c_a^{n_a+1}}{(n_a+1)\big(-\!\!\la)}}
\newcommand{\cb}{\prod_b\fr{c_b^{n_b+1}}{(n_b+1)\big(-\!\!\lb)}}
\newcommand{\cc}{\prod_c\fr{c_c^{n_c+1}}{(n_c+1)\big(-\!\!\lc)}}
\newcommand{\Aa}{\stx{a}{A}{}\!\!}
\newcommand{\Aalp}{\stx{\alpha}{A}{}\!\!}
\newcommand{\olomea}{\stx{a}{\ol\omega}{}\!\!}
\newcommand{\olomeb}{\stx{b}{\ol\omega}{}\!\!}
\newcommand{\olgma}{\stx{a}{\ol\Gamma}{}\!\!\!}
\newcommand{\olgmb}{\stx{b}{\ol\Gamma}{}\!\!\!}

Now let $r,s$ be two nonnegative integers with $K:=r+s\geq 2$ and $x_\alpha:M^{n_\alpha}_\alpha\to\bbr^{n_\alpha+1}$, $1\leq \alpha\leq s$, be hyperbolic affine hyperspheres of dimension $n_\alpha>0$ with affine mean curvatures $\stx{\alpha}{L}\!\!_1$ and with the origin their common affine center. For convenience we make the following convention:
$$1\leq a,b,c\cdots\leq K,\quad 1\leq\lambda,\mu,\nu\leq K-1,\quad
1\leq\alpha,\beta,\gamma\leq s,\quad \td\alpha=\alpha+r,\ \td\beta=\beta+r,\ \td\gamma=\gamma+r.
$$
Furthermore, for each $\alpha=1,\cdots,s$, set $\td i_{\alpha}=i_\alpha+K-1+\sum_{\beta<\alpha}n_\beta$ with $1\leq i_\alpha\leq n_\alpha$.

Define
$$
f_a:=\begin{cases} a,&1\leq a\leq r;\\ \sum_{\beta\leq \alpha}n_\beta+\td{\alpha},&r+1\leq a=\td\alpha\leq r+s,
\end{cases}
$$
and
$$e_a:=\exp\left(-\fr{t_{a-1}}{n_{a}+1}+\fr{t_{a}}{f_{a}}+\fr{t_{a+1}}{f_{a+1}} +\cdots+\fr{t_{K-1}}{f_{K-1}}\right),\quad 1\leq a\leq K=r+s$$
In particular,
$$
e_1=\exp\left(\fr{t_1}{f_1}+\fr{t_2}{f_2} +\cdots+\fr{t_{K-1}}{f_{K-1}}\right),\quad
e_K=\exp\left(-\fr{t_{K-1}}{n_K+1}\right).
$$

Put $n=\sum_\alpha n_\alpha+K-1$ and $M^n=R^{K-1}\times M^{n_1}_1\times\cdots\times M^{n_s}_s$. For any $K$ positive numbers $c_1,\cdots,c_K$, define a smooth map $x:M^n\to\bbr^{n+1}$ by
\begin{align}
x(t^1,&\cdots,t^{K-1},p_1,\cdots,p_s):=(c_1e_1,\cdots, c_re_r,c_{r+1}e_{r+1}x_1(p_1),\cdots,c_Ke_Kx_s(p_s)),\nnm\\&\hs{1cm}\forall (t^1,\cdots,t^{K-1},p_1,\cdots,p_s)\in M^n.\label{mulpro2}
\end{align}

{\prop\label{general sense}\cite{lix11} The map $x:M^n\to\bbr^{n+1}$ defined above is a new hyperbolic affine hypersphere with the affine mean curvature
\be\label{newl1c}
L_1=-\fr1{(n+1)C},\quad C:=\left(\fr1{n+1}\prod_{a=1}^r c_a^2\cdot\prod_{\alpha=1}^s\fr{c_{r+\alpha}^{2(n_\alpha+1)}} {(n_\alpha+1)^{n_\alpha+1}(-\!\!\stx{\alpha}{L}_1)^{n_\alpha+2}}\right)^{\fr1{n+2}},
\ee
Moreover, for given positive numbers $c_1,\cdots,c_K$, there exits some $c>0$ and $c'>0$ such that
the following three hyperbolic affine hyperspheres
\bea &x:=(c_1e_1,\cdots, c_re_r,c_{r+1}e_{r+1}x_1,\cdots,c_Ke_sx_s),\nnm\\
&\bar x:=c(e_1,\cdots, e_r,e_{r+1}x_1,\cdots,e_sx_s),\nnm\\
&\td x:=(e_1,\cdots, e_r,e_{r+1}x_1,\cdots,c'e_sx_s)\nnm
\eea
are equiaffine equivalent to each other.}

{\dfn\label{df2}\cite{lix11} \rm The hyperbolic affine hypersphere $x$ is called the Calabi composition of $r$ points and $s$ hyperbolic affine hyperspheres.}

\rmk\rm The special two cases of the above proposition when $r=0,s=2$ and $r=s=1$, respectively, are discussed in \cite{dil-vra94} and \cite{hu-li-vra08}.

Denote by $\{v^{i_\alpha}_\alpha; \ i_\alpha=1,\cdots,n_\alpha\}$ the local coordinate system of $M_\alpha$, $\alpha=1,\cdots,s$. Then we have

{\prop\label{corr0}\cite{lix11} The Berwald-Blaschke metric $g$, the affine mean curvature $L_1$ and the possibly nonzero components of the Fubini-Pick form $A$ of the Calabi composition $x:M^n\to \bbr^{n+1}$ of $r$ points and $s$ hyperbolic affine hyperspheres $x_\alpha:M_\alpha\to\bbr^{n_\alpha+1}$, $\alpha=1,\cdots,s$, are given as follows:
\begin{align}
&g_{\lambda\mu}=\begin{cases}\displaystyle \fr{\lambda+1}{\lambda}C\delta_{\lambda\mu},&1\leq\lambda\leq r-1;\\
\displaystyle\fr{n_1+r+1}{r}C\delta_{r\mu},&\lambda=r;\\
\displaystyle\fr{\sum_{\beta\leq\alpha+1}n_\beta+\td{\alpha}+1}  {(n_\alpha+1)(\sum_{\beta\leq \alpha}n_\beta+\td{\alpha})}C\delta_{\lambda\mu}, &r+1\leq\lambda=\td\alpha\leq r+s-1.
\end{cases}
\\
&g_{\td i_{\alpha}\td j_{\beta}}=(n_\alpha+1) (-\!\!\lalp)C\galp_{i_{\alpha}j_{\alpha}}\delta_{\alpha\beta},
\quad
g_{\lambda\td i_{\td\alpha}}=0.\label{g-gen}
\\
&A_{\lambda\lambda\lambda} =\begin{cases}\displaystyle \fr{1-\lambda^2}{\lambda^2}C,&1\leq\lambda\leq r-1,\\
\displaystyle\left(\fr1{r^2}-\fr1{(n_1+1)^2}\right)C,&\lambda=r,\\
\displaystyle\fr{(\sum_{\beta\leq\alpha+1}n_\beta+\td{\alpha}+1)C} {(n_{\alpha+1}+1)(\sum_{\beta\leq \alpha}n_\beta+\td{\alpha})}\left(\fr1{\sum_{\beta\leq \alpha}n_\beta+\td{\alpha}}-\fr1{n_{\alpha+1}+1}\right), &r+1\leq\lambda=\td\alpha\leq r+s-1.
\end{cases}
\\
&A_{\lambda\lambda\mu} =\begin{cases}\displaystyle \fr{\lambda+1}{\lambda\mu}C,&1\leq\lambda<\mu\leq r,\\
\displaystyle\fr{(\lambda+1)C}{\lambda(\sum_{\beta\leq \alpha}n_\beta+\td{\alpha})},&1\leq\lambda\leq r-1, \mu=\td\alpha,\\
\displaystyle\fr{(n_1+r+1)C}{r(\sum_{\beta\leq \alpha}n_\beta+\td{\alpha})},&\lambda=r,\ \mu=\td\alpha,\\
\displaystyle\fr{(\sum_{\gamma\leq \alpha+1}n_\gamma+\td{\alpha}+1)C} {(n_{\alpha+1}+1)(\sum_{\gamma\leq \alpha}n_\gamma+\td{\alpha})(\sum_{\gamma\leq \beta}n_\gamma+\td{\beta})},&r+1\leq\lambda=\td\alpha<\mu=\td\beta\leq r+s-1.
\end{cases}
\\
&A_{\td i_{\alpha}\td j_{\alpha}\,{\td\alpha}-1} =-\fr1{n_\alpha+1}g_{\td i_{\alpha}\td j_{\alpha}}=-\!\lalp C\galp_{i_\alpha j_\alpha},\\
&A_{\td i_{\alpha}\td j_{\alpha}\td\beta}=\fr1{\sum_{\gamma\leq \beta}n_\gamma+\td{\beta}}g_{\td i_{\alpha}\td j_{\alpha}}=\fr{(n_\alpha+1)\big(-\!\!\lalp\big)C}{\sum_{\gamma\leq \beta}n_\gamma+\td\beta}\ \galp_{i_\alpha j_\alpha},\quad \beta\geq \alpha,\\
&A_{\td i_{\alpha}\td j_{\alpha}\td k_{\alpha}}=(n_\alpha+1)\big(-\!\!\lalp\big)C\Aalp_{i_\alpha j_\alpha k_\alpha},
\end{align}
where $\stx{\alpha}{L}_1$, $\stx{\alpha}{g}$ and $\stx{\alpha}A$ are the affine mean curvature, the Berwald-Blaschke metric and the Fubini-Pick form of $x_\alpha$, $\alpha=1,\cdots,s$.}

From Proposition \ref{corr0}, the following corollary is easily derived:

{\cor The Calabi composition $x:M^n\to \bbr^{n+1}$ of $r$ points and $s$ hyperbolic affine hyperspheres $x_\alpha:M_\alpha\to\bbr^{n_\alpha+1}$, $\alpha=1,\cdots,s$, is of parallel Fubini-Pick form if and only if for each $\alpha$, the Fubini-Pick form $A_\alpha$ of the factor $x_\alpha$ is parallel.}

By restrictions, $g$ defines a flat metric $g_0$ on $\bbr^{K-1}$ with matrix $(g_{\lambda\mu})$ and, for each $\alpha$, a metric $g_\alpha$ on $M_\alpha$ with matrix $\big(g^\alpha_{i_\alpha j_\alpha}\big)=\big(g_{\td i_{\alpha}\td j_{\alpha}}\big)$ and inverse matrix $\big(g^{i_\alpha j_\alpha}_\alpha\big)$, which is conformal to the original metric $\stx{\alpha}{g}$, or more precisely, \be\label{g_alpha}g_\alpha=(n_\alpha+1) \big(-\!\!\stx{\alpha}{L}_1)C\stx{\alpha}{g}.\ee

\expl\label{expl}\rm Given a positive number $C_0$, let $x_0:\bbr^{n_0}\to \bbr^{n_0+1}$ be the well known flat hyperbolic affine hypersphere of dimension $n_0$ which is defined by
$$
x^1\cdots x^{n_0} x^{n_0+1}=C_0,\quad x^1>0,\cdots,x^{n_0+1}>0.
$$
Then it is not hard to see that $x_0$ is the Calabi composition of $n_0+1$ points. In fact, we can write for example
$$
x_0=(e_1,\cdots,e_{n_0},C_0e_{n_0+1}).
$$
Then by Proposition \ref{corr0} the Berwald-Blaschke metric $g_0$, the affine mean curvature $\stx{0}{L}_1$ and the Fubini-Pick form $\stx{0}{A}$ of $x_0$ are respectively given by (cf. \cite{li-sim-zhao93})
\begin{align}
\stx{0}{g}_{\lambda\mu}=& \fr{\lambda+1}{\lambda}\left(\fr{C_0^2}{n_0+1}\right)^{\fr1{n_0+2}} \delta_{\lambda\mu},\label{gofexpl}\\
\stx{0}{L}_1=&-\fr1{(n_0+1)C} =-(n_0+1)^{-\fr{n_0+1}{n_0+2}}C_0^{-\fr2{n_0+2}},\label{l1ofexpl}\\
\stx{0}{A}_{\lambda\mu\nu}=&\begin{cases}-\fr{\lambda^2-1}{\lambda^2} \left(\fr{C_0^2}{n_0+1}\right)^{\fr1{n_0+2}}, &{\rm if\ }\lambda=\mu=\nu;\\
\fr{\lambda+1}{\lambda\nu}\left(\fr{C_0^2}{n_0+1}\right)^{\fr1{n_0+2}}, &{\rm if\ }\lambda=\mu<\nu;\\
0, &{\rm otherwise.}
\end{cases}\label{aofexpl}
\end{align}
Thus the Pick invariant of $x_0$ is
\be\label{jofexpl}
\stx{0}{J}=\fr1{n_0(n_0-1)}\stx{0}{g}{}\!\!^{\lambda_1\lambda_2} \stx{0}{g}{}\!\!^{\mu_1\mu_2} \stx{0}{g}{}\!\!^{\nu_1\nu_2} \stx{0}{A}_{\lambda_1\mu_1\nu_1}\stx{0}{A}_{\lambda_2\mu_2\nu_2} =(n_0+1)^{-\fr{n_0+1}{n_0+2}}C_0^{-\fr2{n_0+2}}=-\!\!\stx{0}{L}_1.
\ee

To end this section we list some properties of the Calabi composition of points and hyperbolic affine hyperspheres.

Write $M_0=\bbr^{K-1}$. Then, with respect to the metric $g$ on $M^n$, the Fubini-Pick form $A$ can be identified with a $TM^n$-valued symmetric $2$-form $\sigma:TM^n\times TM^n\to TM^n$. For each ordered triple $\alpha,\beta,\gamma\in\{0,1,\cdots,s\}$, $\sigma$ defines one $TM_\gamma$-valued bilinear map $\sigma^{\gamma}_{\alpha\beta}:TM_{\alpha}\times TM_{\beta}\to TM_\gamma$, which is the $TM_\gamma$-component of $\sigma_{\alpha\beta}$, the restriction of $\sigma$ to $TM_{\alpha}\times TM_{\beta}$. Define
$$
H_\alpha=\fr1{n_{\alpha}}\tr_{g_\alpha}\sigma^0_{\alpha\alpha}\equiv \fr1{n_{\alpha}}g^{i_\alpha j_\alpha}_\alpha\sigma^0_{\alpha\alpha}\left(\pp{}{v^{i_{\alpha}}_{\alpha}}, \pp{}{v^{j_{\alpha}}_{\alpha}}\right),
$$
where the metric $g_\alpha$ is given by \eqref{g_alpha}.

{\prop\label{properties}\cite{lix11} Let $x:M^n\to \bbr^{n+1}$ be the Calabi composition of $r$ points and $s$ hyperbolic affine hyperspheres and $g$ the Berwald-Blaschke metric of $x$. Then

$(1)$ The Riemannian manifold $M^n\equiv (M^n,g)$ is reducible;

$(2)$ There must be a positive dimensional Euclidean factor $\bbr^q$ in the de Rham decomposition of $M^n$;

$(3)$ $q\geq s-1$ with the equality holding if and only if $r=0$;

$(4)$ $\sigma^\gamma_{\alpha\beta}\equiv 0$ if
$(\alpha,\beta,\gamma)$ is not one of the following triples: $(0,0,0)$, $(\alpha,\alpha,0)$, $(\alpha,0,\alpha)$, $(0,\alpha,\alpha)$ or $(\alpha,\alpha,\alpha)$.

$(5)$ The vector-valued functions $H_\alpha$, $\alpha=1,\cdots,s$, satisfies the following qualities:
\begin{align}
&g(H_\alpha,H_\alpha)=C^{-1}\left(\fr1{n_{\alpha}+1}-\fr1{f_K}\right) =\fr{n-n_\alpha}{n_\alpha+1}(-L_1),\\ &g(H_\alpha,H_\beta)=L_1\quad{\rm for\ }\alpha\neq\beta;
\end{align}

$(6)$ $\sigma^{\alpha}_{\alpha\alpha}$ is identical to the $TM_\alpha$-valued symmetric bilinear form defined by the Fubini-Pick form $\stx{\alpha}{A}$ of $x_\alpha$.}

\rmk\rm In the next section, we shall show that a locally strongly convex hypersurface $x:M^n\to \bbr^{n+1}$ with parallel Fubini-Pick form is locally the Calabi composition of some points and hyperbolic affine hyperspheres if and only if the above condition (1) holds (see Theorem \ref{chara}).

\section{A duality correspondence theorem}

\subsection{The correspondence theorem}

In this section, we prove a theorem which locally establishes a simple correspondence between the set of symmetric affine hypersurfaces and that of the minimal symmetric Lagrangian submanifolds immersed in some complex space form. This is one of the key results that may provide another way to establish the classification for those important hypersurfaces.

{\thm\label{corr thm} Let $x:M^n\to \bbr^{n+1}$ be a simply connected and locally strongly convex symmetric equiaffine hypersphere of affine mean curvature $L_1$. Then $x$ defines, uniquely up to certain equivalences, a simply connected, minimal symmetric Lagrangian submanifold $\td x: \td M^n\to \bbq^n(-4L_1)$ immersed in the complex space form $\bbq^n(-4L_1)$ of constant holomorphic sectional curvature $-4L_1$; Conversely, each of the simply connected, minimal symmetric Lagrangian submanifolds $\td x$ immersed in the complex space form $\bbq^n(4c)$  of constant holomorphic sectional curvature $4c$ corresponds, uniquely up to affine equivalences, to one simply connected and locally strongly convex symmetric equiaffine hypersphere $x$ of affine mean curvature $-c$, which defines $\td x$.

In other words, there is a one to one correspondence between the set of the equiaffine equivalence classes of simply connected and locally strongly convex symmetric affine hyperspheres of affine mean curvature $L_1$ and that of the equivalence classes of simply connected, minimal symmetric Lagrangian submanifolds immersed in the complex space form $\bbq^n(-4L_1)$ of constant holomorphic sectional curvature $-4L_1$.}

\proof The assumption that $x$ is an equiaffine hypersphere implies that the affine fundamental form $B$ is a scalar multiple of the Berwald-Blaschke metric $g$, that is, $B=L_1\,g$. Then \eqref{gaus} and \eqref{afsp} show that the affine Riemannian curvature tensor can be rewritten as
\be\label{gaus1}
R(X,Y)Z=L_1(g(Y,Z)X-g(X,Z)Y)-[A(X),A(Y)]Z.
\ee
On the other hand, since $x$ is symmetric, we can write $M^n=G/K$. As the symmetric space, $M^n$ has a unique dual space $\td M^n=\td G/K$ which is also simply connected.

By the definition of the dual space of a symmetric space, the Lie algebras $\mathfrak g$, $\td{\mathfrak g}$ of $G$, $\td G$ can be decomposed respectively into
$$
{\mathfrak g}={\mathfrak k}+{\mathfrak m}, \quad \td{\mathfrak g}={\mathfrak k}+\sqrt{-1}{\mathfrak m},
$$
where ${\mathfrak k}$ is the Lie algebra of $K$, and ${\mathfrak m}$, $\td{\mathfrak m}:=\sqrt{-1}{\mathfrak m}$ are respectively identified with the tangent spaces $T_oM$, $T_{\td o}\td M^n$ at the base points $o:=eK$, $\td o:=\td eK$, where $e,\td e$ are the unit elements of $G$ and $\td G$, respectively.

Clearly, since the Fubini-Pick form $A$ of the symmetric affine hypersphere $x$ is $G$-invariant, by taking the value at the point $o$, it defines a $K$-invariant symmetric trilinear form $A_o$ on ${\mathfrak m}$ which is identified with a $K$-invariant bilinear map $A_o:{\mathfrak m}\times{\mathfrak m}\to {\mathfrak m}$ by using the metric $g$ at $o$. More precisely
$$
\lagl A_o(X,Y),Z\ragl=A_o(X,Y,Z),\quad\forall\,X,Y,Z\in {\mathfrak m}.
$$
The invariance of $A_o$ by $K$ is equivalent to that ${\mathfrak k}\cdot A_o=0$. As mentioned earlier, we also take $A_o$ as a map of ${\mathfrak m}$ into $\ed({\mathfrak m})$, that is, for each $X\in {\mathfrak m}$, we have a linear map $A_X:{\mathfrak m}\to {\mathfrak m}$ given by $A_X(Y):=A_o(X,Y)$.

Taking the linear isomorphism $\sqrt{-1}:{\mathfrak m}\to \sqrt{-1}{\mathfrak m}$ to be isometric, one defines a $K$-invariant scalar inner product on $T_{\td o}\td M^n$ which in turn determines a $\td G$-invariant metric $\td g$ on $\td M^n$. This invariant metric has a curvature tensor $\td R$ and, at $\td o$, it is the minus of $R$ at $o$ under the identification map $\sqrt{-1}$.

On the other hand, by using the identification $\sqrt{-1}:{\mathfrak m}\to \sqrt{-1}{\mathfrak m}$, we can define a $K$-invariant symmetric trilinear form $\td\sigma$ on $\td{\mathfrak m}$ by
\be\label{sigmatri}
\hs{-.3cm}\td\sigma(\sqrt{-1}X,\sqrt{-1}Y,\sqrt{-1}Z)=A_o(X,Y,Z),\,
\forall\,X,Y,Z\in{\mathfrak m}.
\ee

The corresponding symmetric bilinear map $\td \sigma:\td{\mathfrak m}\times \td{\mathfrak m}\to \td{\mathfrak m}$ is given by
\be\label{sigmabi}
\hs{-1cm}\td\sigma(\sqrt{-1}X,\sqrt{-1}Y)=\sqrt{-1}A_o(X,Y),\quad
\forall\,X,Y\in{\mathfrak m}.
\ee
Then it is not hard to see that, for any $\td X,\td Y,\td Z\in\td{\mathfrak m}$,
\be\label{gaus2}
\td R(\td X,\td Y)\td Z=(-L_1)(g(\td Y,\td Z)\td X-g(\td X,\td Z)\td Y)+[\td\sigma_{\td X},\td\sigma_{\td Y}](\td Z).
\ee
Furthermore, the $K$-invariance of $\td\sigma$ is equivalent to that ${\mathfrak k}\cdot\td\sigma=0$. From this we see that $\td\sigma$ uniquely defines a totally symmetric $\td G$-invariant bilinear form $\td\sigma\in\Gamma\big((\bigodot^2 T^*\td M^n)\otimes T\td M^n)$. Thus the equation \eqref{gaus2} holds globally on $\td M^n$.

Now by the same argument of \cite{nai81} we know that, up to holomorphic isometries on $\bbq^n(-4L_1)$, there is a unique minimal symmetric Lagrangian submanifold $\td x:\td M^n\to \bbq^n(-4L_1)$, such that the above $\td\sigma$ coincides with the one induced by the second fundamental form of $\td x$.

Conversely, let $\td x:\td M^n\to \bbq^n(4c)$ be a simply connected, minimal symmetric Lagrangian submanifold with metric $\td g$. Suppose that the dual space of $\td M^n$ is $M^n$, on which a unique invariant metric $g$ is naturally determined by $\td g$. Then we can write $\td M^n=\td G/K$ and $M^n=G/K$ for suitable Lie groups $\td G$ and $G$ where $K$ is one common closed Lie subgroup of $\td G$ and $G$. Write $\td{\mathfrak m}=T_{\td o}\td M^n$ with $\td o=\td eK\in \td M^n$. Then ${\mathfrak m}:=\sqrt{-1}\td{\mathfrak m}\equiv T_oM$ with $o=eK\in M^n$.

It is not hard to see that the second fundamental form $\sigma$ of $\td x$ at $\td o$ uniquely defines a totally symmetric trilinear form $\td\sigma:\td{\mathfrak m}\times\td{\mathfrak m}\times\td{\mathfrak m}\to\bbr$ which gives a totally symmetric trilinear form $A_o:{\mathfrak m}\times{\mathfrak m}\times{\mathfrak m}\to\bbr$, identified with one $\mathfrak{m}$-valued symmetric $2$-form $A_o:\mathfrak{m}\times\mathfrak{m}\to\mathfrak{m}$ and one linear map $A_X:\mathfrak{m}\to \mathfrak{m}$ for each $X\in\mathfrak{m}$.
The curvature tensor $R$ of $M^n$ at $o$ is the minus of the
curvature $\td R$ of $\td M^n$ at $\td o$. It then follows from the Gaussian equation of
$\td x$ (see \eqref{gauss3}) that
\be\label{gaus3}
R(X,Y)Z=-c(g(Y,Z)X-g(X,Z)Y)-[A_X,A_Y]Z,
\ee
for all $X,Y,Z\in {\mathfrak m}$.

Since $\td\sigma$ is $\td G$-invariant we have ${\mathfrak k} \cdot\td\sigma=0$ which is equivalent to that ${\mathfrak k} \cdot A_o=0$. This implies that the trilinear form $A_o$ given by $\td\sigma$ extends to a globally defined, $G$-invariant and totally symmetric trilinear form $A$. Therefore \eqref{gaus2} holds everywhere on $M^n$ since $R$ is also $G$-invariant.

Now we can apply the existence and uniqueness theorem of equiaffine geometry of hypersurfaces (Theorems \ref{affine existence} and \ref{affine uniqueness}) to conclude that, up to affine equivalences, there exists uniquely one locally strongly convex symmetric affine hypersphere $x:M^n\to \bbr^{n+1}$ of which the Berwald-Blaschke metric and the Fubini-Pick form coincide with
the above $G$-invariant metric $g$ and the trilinear form $A$.
\endproof

Similar to ${\mathcal S}_{M^n}(c)$ given in Section 2, we define
$\td{\mathcal S}_{M^n}(L_1)$ via the basic equations \eqref{gaus_af sph}, \eqref{basic1} and the condition ${\mathfrak k}\cdot A=0$.
Then the following corollary is directly derived from Theorem \ref{corr thm} and Proposition \ref{nai81}:

{\cor\label{corr1} Let $(M^n,g)= (G/K,g)$ be a simply connected Riemannian symmetric space of dimension $n$ with  $\td{\mathcal S}_{M^n}(L_1)\neq\emptyset$ with the symmetric metric $g$. Then for each $A\in \td {\mathcal S}_{M^n}(L_1)$, there exists uniquely one locally strongly convex symmetric affine hypersphere $x:M^n\to \bbr^{n+1}$, such that the corresponding Berwald-Blaschke metric coincides with $g$ and the Fubini-Pick form is given by the $T_oM^n$-valued symmetric bilinear form $A$.}

\subsection{A characterization of Calabi product}

A direct application of Theorem \ref{corr thm} is to establish a necessary and sufficient condition for a locally strongly convex hypersurface with parallel Fubini-Pick form locally to be the Calabi composition of several hyperbolic affine hyperspheres, possibly including point factors. It turns out that
this special characterization theorem is also needed in the next section for proving the main classification theorem. Here we should remark that, in \cite{hu-li-vra08}, Z.J. Hu etc give a characterization of the two factor Calabi composition but in a different manner.

Note that for a given locally strongly convex hypersurface $x:M^n\to\bbr^{n+1}$ with the Berwald-Blaschke metric $g$, $(M^n,g)$ is a Riemannian manifold.

{\thm\label{chara} A locally strongly convex hypersurface $x:M^n\to\bbr^{n+1}$ with parallel Fubini-Pick form is locally affine equivalent to the Calabi composition of some hyperbolic affine hyperspheres possibly including point factors if and only if $M^n$ is reducible as a Riemannian manifold with respect to the Berwald-Blaschke metric.}

\proof
The necessary part of the theorem is obvious (cf. Proposition \ref{properties}). To prove the sufficient part, we first use a result in \cite{dil-vra-yap94} to know that, under the assumption of the theorem, $x$ must be a hyperbolic affine hypersphere. Therefore, the affine mean curvature $L_1$ is a negative constant. Note that, by Proposition \ref{sym para1}, $x$ is locally symmetric as an equiaffine hypersphere. Therefore, without loss of generality, we can assume that $x$ is a simply connected symmetric equiaffine hypersphere. Thus, by Theorem \ref{corr thm}, $x$ uniquely defines a minimal symmetric Lagrangian submanifold $\td x:\td M\to \bbc P^n(-4L_1)$ immersed in the complex projective space $\bbc P^n(-4L_1)$ with constant holomorphic sectional curvature $-4L_1$, where $\td M$ is the dual symmetric space of $M^n$ and the second fundamental form $\td\sigma$ is determined by the Fubini-Pick form $A$. Since $M^n$ is reducible, $\td M$ is also reducible. By Propostion \ref{sym para2}, $\td x$ is parallel as an immersion. It then follows from Lemma 4.1 in \cite{nai81} that $\td M$ must have an Euclidean factor $\bbr^{n_0}$, $n_0>0$, in its de Rham decomposition:
$$
\td M=\bbr^{n_0}\times\td M_1\times\cdots\times\td M_s,
$$
where $M_1,\cdots M_s$ are simply connected compact symmetric spaces. Thus, if we write $\td M=\td G/K$, then the Lie algebras $\td{\mathfrak g},{\mathfrak k}$ of $\td G$ and $K$ have respectively the following decompositions
\bea
&\td{\mathfrak g}={\mathfrak g_0}\oplus \td{\mathfrak g}_1\oplus\cdots \td{\mathfrak g}_s={\mathfrak k}\oplus \td{\mathfrak m},\\
&{\mathfrak k}={\mathfrak k}_0\oplus {\mathfrak k}_1\oplus\cdots {\mathfrak k}_s,
\eea
where $({\mathfrak g}_\alpha,{\mathfrak k}_\alpha)$ is the symmetric pair of Lie algebras corresponding to the symmetric factor $\td M_\alpha$ for $\alpha=1,\cdots,s$ and ${\mathfrak m}_0:=T_o\bbr^{n_0}\equiv \bbr^{n_0}$, ${\mathfrak m}=T_{\td o}\td M$ with $o$ the origin of $\bbr^{n_0}$ and $\td o=\td eK$. Therefore for each $\alpha$, $\td{\mathfrak g}_\alpha$ is decomposed into
$\td{\mathfrak g}_\alpha={\mathfrak k}_\alpha\oplus {\mathfrak m}_s$ with ${\mathfrak m}_\alpha=T_o\td M_\alpha$. It follows that ${\mathfrak m}={\mathfrak m}_0\oplus {\mathfrak m}_1\oplus\cdots\oplus {\mathfrak m}_s$. The second fundamental form $\td\sigma$ defines a $TM^*$-valued symmetric bilinear form, still denoted by $\td\sigma$, of which the restriction to the given point $\td o$ is a $\td{\mathfrak m}$-valued symmetric bilinear form on $\td{\mathfrak m}$.
By Theorem 6.4 in \cite{nai81}, $\td\sigma$ can be decomposed as
\be\label{decomposition tdsigma}
\td\sigma=\sum_{\alpha=0}^s\td\sigma_{\alpha\alpha}^\alpha +\sum_{\alpha=1}^s\td\sigma_{\alpha\alpha}^0 +\sum_{\alpha=1}^s\td\sigma_{\alpha 0}^\alpha+\sum_{\alpha=1}^s\td\sigma_{0\alpha}^\alpha
\ee
where, for each triple $(\alpha,\beta,\gamma)$, the bilinear map $\td\sigma_{\alpha\beta}^\gamma:\td{\mathfrak m}_\alpha\times\td{\mathfrak m}_\beta\to\td{\mathfrak m}_\gamma$ is the $\td{\mathfrak m}_\gamma$-component of $\td\sigma$ restricting to the subspace $\td{\mathfrak m}_\alpha\times\td{\mathfrak m}_\beta$.

Following Naitoh \cite{nai81}, {\em the $\alpha$-th mean curvature} $H_\alpha$ is defined to be the $\fr1{n_\alpha}$ multiple of the trace of $\td\sigma^0_{\alpha\alpha}$ with respect to the metric $\td g_\alpha$ on $\td M_\alpha$, that is, $H_\alpha=\fr1{n_\alpha}\tr_{\td g_\alpha}(\td\sigma^0_{\alpha\alpha})$. Denote $\td c_\alpha=| H_\alpha|$. Then by Naitoh (Theorem 6.4,\cite{nai81}) together with \eqref{decomposition tdsigma} we have the following conclusions:
\begin{align}
&\td\sigma^\alpha_{\alpha\alpha}\in {\mathcal S}_{\td M_\alpha}(-L_1+\td c^2_\alpha),\label{alp alp alp}\\
&\td\sigma^0_{00}\in \ol{\mathcal S}_{\bbr^{n_0}}(-L_1),\quad\td\sigma^0_{00}(Z_0,H_\alpha) =g(Z_0,H_\alpha)H_\alpha+L_1Z_0,\label{000}\\
&\td\sigma^0_{\alpha\alpha}(X_\alpha,Y_\alpha)=g(X_\alpha,Y_\alpha)H_\alpha,\quad g(H_\alpha,H_\beta)=L_1,\quad{\rm if\ }1\leq\alpha\neq\beta\leq s,\label{alp alp 0}\\
&\td\sigma^\alpha_{\alpha 0}(X_\alpha,Z_0) =\td\sigma^\alpha_{0\alpha}(Z_0,X_\alpha)=g(Z_0,H_\alpha)X_\alpha.\label{alp 0 alp}
\end{align}

Define
${\mathcal H}=\spn\{H_1,\cdots H_s\}$ and denote by
${\mathcal H}^\bot$ the orthogonal complement of ${\mathcal H}$ in ${\mathfrak m}_0$. Thus ${\mathfrak m}_0={\mathcal H}\oplus {\mathcal H}^\bot$, and $\td\sigma^0_{00}$ can be decomposed into the sum of its ${\mathcal H}$-component $\td\sigma_0^{\mathcal H}$ and its ${\mathcal H^\bot}$-component $\td\sigma_0^{\mathcal H^\bot}$, that is,
$\td\sigma^0_{00}=\td\sigma_0^{\mathcal H}+\td\sigma_0^{\mathcal H^\bot}$.

{\lem\label{rh geq s-1} Let $n_0$, $s$ be as above. Then
$n_0\geq \dim {\mathcal H}\geq s-1$. Furthermore, $n_0\geq s$ if and only if $\dim {\mathcal H}=s$.}

{\it Proof of Lemma \ref{rh geq s-1}}:

To prove the first part of the lemma, it suffices to show that the set of the $s$ nonzero vectors $H_1,\cdots,H_s$ in ${\mathfrak m}_0\equiv T_o\bbr^{n_0}$ has a rank not less than $s-1$. This is equivalent to show that the $s$-th order matrix
$$
\lmx g(H_1,H_1)&g(H_1,H_2)&\cdots&g(H_1,H_s)\\
g(H_2,H_1)&g(H_2,H_2)&\cdots&g(H_2,H_s)\\
\cdots&\cdots&\cdots&\cdots\\
g(H_s,H_1)&g(H_s,H_2)&\cdots&g(H_s,H_s)
\rmx
=\lmx \td c^2_1&L_1&\cdots&L_1\\
L_1&\td c^2_2&\cdots&L_1\\
\cdots&\cdots&\cdots&\cdots\\
L_1&L_1&\cdots&\td c^2_s
\rmx
$$
has a rank equal or larger than $s-1$. Indeed, by deleting the last line and the second last column, the above matrix has a ($s-1$)-minor
\begin{align}
&\det\lmx \td c^2_1&L_1&\cdots&L_1\\
L_1&\td c^2_2&\cdots&L_1\\
\cdots&\cdots&\cdots&\cdots\\
L_1&L_1&\cdots&L_1
\rmx
=L_1\det\lmx \td c^2_1&L_1&\cdots&L_1\\
L_1&\td c^2_2&\cdots&L_1\\
\cdots&\cdots&\cdots&\cdots\\
1&1&\cdots&1
\rmx\nnm\\
=&L_1\det\lmx \td c^2_1-L_1&0&\cdots&0\\
0&\td c^2_2-L_1&\cdots&0\\
\cdots&\cdots&\cdots&\cdots\\
1&1&\cdots&1
\rmx
=L_1(\td c^2_1-L_1)\cdots(\td c^2_{s-2}-L_1)<0,\nnm
\end{align}
where we have used the fact that $L_1<0$.

Furthermore, if $n_0\geq s$, and $\dim {\mathcal H}=s-1$, then $r-1:=\dim{\mathcal H}^\bot\geq 1$. Consider the restriction $\bar \sigma_0^{\mathcal H}$ of $\td \sigma_0^{\mathcal H}$ to the subspace ${\mathcal H}^\bot\times{\mathcal H}^\bot$. Define $H_0=\fr1{r-1}\tr\bar \sigma_0^{\mathcal H}$ and $\td c_0=|H_0|$. Then, for any unit vector $e_0\in {\mathcal H}^\bot$ and each $\alpha=1,\cdots,s$, we have
$$
g(\bar\sigma^{\mathcal H}_0(e_0,e_0),H_\alpha)=g(\td\sigma(e_0,e_0),H_\alpha) =g(\td\sigma^0_{00}(e_0,H_\alpha),e_0)=L_1g(e_0,e_0)=L_1,
$$
implying that
$$
g(H_0,H_\alpha)=L_1,\quad \alpha=1,\cdots,s.
$$
Then in the same way as in proving that the rank of the matrix $(g(H_\alpha,H_\beta))_{1\leq\alpha,\beta\leq s}$ of order $s$ is no less than $s-1$ we can obtain that the rank of the $(s+1)$-th order matrix $(g(H_\alpha,H_\beta))_{0\leq\alpha,\beta\leq s}$ is no less than $s$. Since $\{H_\alpha;\ 0\leq\alpha\leq s\}\subset {\mathcal H}$, it follows that $\dim {\mathcal H}\geq s$ which contradicts the assumption.\endproof

Define $\bar\sigma_0^{\mathcal H^\bot}=\left.\td\sigma_0^{\mathcal H^\bot}\right|_{{\mathcal H^\bot}\times {\mathcal H^\bot}}$. By \eqref{000}, for any $X_{\mathcal H},Y_{\mathcal H}\in {\mathcal H}$, $\td\sigma^0_{00}(X_{\mathcal H},Y_{\mathcal H})\in {\mathcal H}$, and for any $Y_{\mathcal H^\bot}\in {\mathcal H^\bot}$, $\td\sigma^0_{00}(X_{\mathcal H},Y_{\mathcal H^\bot})\in {\mathcal H^\bot}$. Therefore, $\td\sigma^0_{00}$ can be decomposed into the following components:
\be\label{decomp tdsigma0}
\td\sigma^0_{00}=\bar\sigma_0^{\mathcal H^\bot}+\bar\sigma_0^{\mathcal H}+\left.\td\sigma_0^{\mathcal H^\bot}\right|_{{\mathcal H^\bot}\times {\mathcal H}}+\left.\td\sigma_0^{\mathcal H^\bot}\right|_{{\mathcal H}\times {\mathcal H^\bot}}+\left.\td\sigma_0^{\mathcal H}\right|_{{\mathcal H}\times {\mathcal H}}.
\ee

{\lem\label{barsigma0} Define $r=\dim{\mathcal H}^\bot+1$. Then $\bar\sigma_0^{\mathcal H^\bot}\in {\mathcal S}_{\bbr^{r-1}}(-\fr{(n+1)L_1}r)$ if $\dim{\mathcal H}^\bot\geq 1$.}

{\it Proof of Lemma \ref{barsigma0}}:

Since $\dim{\mathcal H}^\bot\geq 1$, it holds by Lemma \ref{rh geq s-1} that
$$
n_0=\dim{\mathcal H}+\dim{\mathcal H}^\bot\geq s-1+1=s.
$$
Making use of Lemma \ref{rh geq s-1} once again we have that $\dim{\mathcal H}=s$ and therefore $\{H_1,\cdots,H_s\}$ is a basis for the linear space ${\mathcal H}$. Set $h_{\alpha\beta}=g(H_\alpha,H_\beta)$, $1\leq\alpha,\beta\leq s$, and $(h^{\alpha\beta})=(h_{\alpha\beta})^{-1}$.

We first compute $\bar\sigma^{\mathcal H}_0$. Write
$$
\bar\sigma^{\mathcal H}_0(X,Y)=\sum C^\alpha_{XY}H_\alpha,\quad \forall X,Y\in {\mathcal H^\bot}.
$$
Then we have
\begin{align}
g(\bar\sigma^{\mathcal H}_0(X,Y),H_\alpha) =&\sum C^\beta_{XY}g(H_\beta,H_\alpha)=\sum C^\beta_{XY}h_{\beta\alpha};\nnm\\
g(\bar\sigma^{\mathcal H}_0(X,Y),H_\alpha) =&g(\td\sigma^0_{00}(X,Y),H_\alpha) =g(\td\sigma^0_{00}(X,H_\beta),Y)\nnm\\
=&g(g(X,H_\alpha)H_\alpha+L_1X,Y)=L_1g(X,Y),\nnm
\end{align}
implying that
$
\sum C^\beta_{XY}h_{\beta\alpha}=L_1g(X,Y)
$
or equivalently $C^\alpha_{XY}=L_1g(X,Y)\sum_{\alpha,\beta} h^{\alpha\beta}$. It follows that
\be\label{barsigma0h}
\bar\sigma^{\mathcal H}_0(X,Y)=L_1g(X,Y)\sum_\beta h^{\alpha\beta}H_\alpha,
\quad\forall X,Y\in {\mathcal H^\bot}.
\ee
Thus we have
\be\label{H0}
H_0=\fr1{r-1}\tr(\bar\sigma^{\mathcal H}_0)=L_1\sum_{\alpha,\beta} h^{\alpha\beta}H_\alpha.
\ee

On the other hand, by using \eqref{decomposition tdsigma} and the fact that $\tr(\td\sigma)=0$, it is seen that
$$
\tr\td\sigma^0_{00}+\sum_\alpha\tr\td\sigma^0_{\alpha\alpha}=0,
$$
which with the decomposition \eqref{decomp tdsigma0} gives
\begin{align}
&\tr(\bar\sigma^{\mathcal H^\bot}_0)=0,\label{tr barsigma0hbot=0}\\
&(r-1)H_0+\sum_{\alpha,\beta}h^{\alpha\beta}\td\sigma^{\mathcal H}_0(H_\alpha,H_\beta) +\sum_\alpha n_\alpha H_\alpha=0.\label{sum Halpha=0}
\end{align}
But by \eqref{000},
$$
\td\sigma^{\mathcal H}_0(H_\alpha,H_\beta)=g(H_\alpha,H_\beta)H_\beta +L_1H_\alpha =h_{\alpha\beta} H_\beta +L_1H_\alpha.
$$
Thus \eqref{sum Halpha=0} can be rewritten as
\be\label{sum Halpha=0-1}
(r-1)H_0+\sum_\alpha(1+L_1\sum_\beta h^{\alpha\beta})H_\alpha +\sum_\alpha n_\alpha H_\alpha=0.
\ee

Comparing \eqref{H0} and \eqref{sum Halpha=0-1} gives that
$$
\sum_\alpha\left(rL_1\sum_\beta h^{\alpha\beta}+(n_\alpha+1)\right)H_\alpha=0,
$$
or equivalently
\be\label{sumbeta h^alphabeta}
rL_1\sum_\beta h^{\alpha\beta}+(n_\alpha+1)=0,\quad \alpha=1,\cdots,s.
\ee
It follows that
\be\label{*}
\sum_\beta h^{\alpha\beta}=-\fr{n_\alpha+1}{rL_1},\quad \forall \alpha
\ee
which with \eqref{barsigma0h} gives
\be\label{barsigma0h-1}
\bar\sigma^{\mathcal H}_0(X,Y)=-\left(\sum_\alpha \fr{n_\alpha+1}{r}\right)g(X,Y)H_\alpha,
\quad\forall X,Y\in {\mathcal H^\bot},
\ee
and thus
\be\label{H0-1}
H_0=\fr1{r-1}\tr(\bar\sigma^{\mathcal H}_0)=-\sum_\alpha \fr{n_\alpha+1}{r}H_\alpha.
\ee

Since we have shown that $\tr(\bar\sigma^{\mathcal H^\bot}_0)=0$ (Equation \eqref{tr barsigma0hbot=0}), to complete the proof of Lemma \ref{barsigma0}, it now suffices to show that
\be\label{barsigma0hbot-gauss}
-\fr{(n+1)L_1}r(g(Y,Z)X-g(X,Z)Y) +[\bar\sigma^{\mathcal H^\bot}_0(X),\bar\sigma^{\mathcal H^\bot}_0(Y)](Z)=0, \quad\forall X,Y,Z\in {\mathcal H^\bot}.
\ee

In fact, for any $X,Y,Z\in {\mathcal H^\bot}$, that $\td\sigma^0_{00}\in \ol{\mathcal S}_{\bbr^{n_0}}(-L_1)$ implies
\be\label{td sigma0 gauss}
-L_1(g(Y,Z)X-g(X,Z)Y)+[\td\sigma^0_{00}(X),\td\sigma^0_{00}(Y)](Z)=0.
\ee
But by the decomposition \eqref{decomp tdsigma0}
\begin{align}
\td\sigma^0_{00}(X)(\td\sigma^0_{00}(Y)(Z)) =&\td\sigma^0_{00}(X,\td\sigma^0_{00}(Y,Z))\nnm\\
=&\td\sigma^0_{00}(X,\bar\sigma^{\mathcal H^\bot}_0(Y,Z) +\bar\sigma^{\mathcal H}_0(Y,Z))\nnm\\
=&\td\sigma^0_{00}(X,\bar\sigma^{\mathcal H^\bot}_0(Y,Z)) +\td\sigma^0_{00}(X,\bar\sigma^{\mathcal H}_0(Y,Z))\nnm\\
=&\bar\sigma^{\mathcal H^\bot}_0(X,\bar\sigma^{\mathcal H^\bot}_0(Y,Z)) +\bar\sigma^{\mathcal H}_0(X,\bar\sigma^{\mathcal H^\bot}_0(Y,Z))
+\td\sigma^{\mathcal H^\bot}_0(X,\bar\sigma^{\mathcal H}_0(Y,Z))\nnm\\
=&\bar\sigma^{\mathcal H^\bot}_0(X)(\bar\sigma^{\mathcal H^\bot}_0(Y)(Z)) -\fr1rg(X,\bar\sigma^{\mathcal H^\bot}_0(Y,Z))\sum_\alpha (n_\alpha+1)H_\alpha\nnm\\ &\quad-\fr1rg(Y,Z)\sum_\alpha (n_\alpha+1)\td\sigma^{\mathcal H^\bot}_0(X,H_\alpha)\nnm\\
=&\bar\sigma^{\mathcal H^\bot}_0(X)(\bar\sigma^{\mathcal H^\bot}_0(Y)(Z))
-\fr1rg(X,\bar\sigma^{\mathcal H^\bot}_0(Y,Z))\sum_\alpha (n_\alpha+1)H_\alpha\nnm\\ &\quad
-\fr1r\sum_\alpha (n_\alpha+1)L_1g(Y,Z)X\nnm\\
=&\bar\sigma^{\mathcal H^\bot}_0(X)(\bar\sigma^{\mathcal H^\bot}_0(Y)(Z))
-\fr1rg(X,\bar\sigma^{\mathcal H^\bot}_0(Y,Z))\sum_\alpha (n_\alpha+1)H_\alpha\nnm\\ &\quad
-\fr{n-r+1}rL_1g(Y,Z)X\nnm
\end{align}
where we have used \eqref{000}, \eqref{barsigma0h-1} and the definition of $r$. Since $g(X,\bar\sigma^{\mathcal H^\bot}_0(Y,Z))=g(Y,\bar\sigma^{\mathcal H^\bot}_0(X,Z))$, we find that
\begin{align}
[\td\sigma^0_{00}(X),\td\sigma^0_{00}(Y)](Z) =&\td\sigma^0_{00}(X)(\td\sigma^0_{00}(Y)(Z)) -\td\sigma^0_{00}(Y)(\td\sigma^0_{00}(X)(Z))\nnm\\
=&\bar\sigma^{\mathcal H^\bot}_0(X)(\bar\sigma^{\mathcal H^\bot}_0(Y)(Z))
-\bar\sigma^{\mathcal H^\bot}_0(Y)(\bar\sigma^{\mathcal H^\bot}_0(X)(Z))\nnm\\
&\ -\fr{n-r+1}rL_1(g(Y,Z)X-g(X,Z)Y)\nnm\\
=&[\bar\sigma^{\mathcal H^\bot}_0(X),\bar\sigma^{\mathcal H^\bot}_0(Y)](Z)
-\fr{n-r+1}rL_1(g(Y,Z)X-g(X,Z)Y).
\end{align}
Inserting the above equality into \eqref{td sigma0 gauss} we obtain the equation \eqref{barsigma0hbot-gauss}, which completes the proof of Lemma \ref{barsigma0}.\endproof

{\lem\label{lem4.3} The vector-valued symmetric bilinear form $\td\sigma\in {\mathcal S}_{\td M^n}(-L_1)$ is uniquely, up to equivalence, determined by the metrics $\td g_\alpha$, the flat metric $g$ on $\bbr^{n_0}$, the bilinear forms $\td\sigma^{\alpha}_{\alpha\alpha}$ $($$\alpha=1,\cdots,s$$)$ and the affine mean curvature $L_1$.}

{\it Proof of Lemma \ref{lem4.3}}:

Since $\td\sigma^{\alpha}_{\alpha\alpha}\in {\mathcal S}_{(\td M_\alpha,\td g_\alpha)}(-L_1+\td c^2_\alpha)$, we see that $\td c_\alpha$ is completely determined by $\td g_\alpha$, $L_1$ and $\td\sigma^\alpha_{\alpha\alpha}$ via
$$
R_{\td g_\alpha}(X_\alpha,Y_\alpha)Z_\alpha=(-L_1+\td c^2_\alpha)(\td g_\alpha(Y_\alpha,Z_\alpha)X_\alpha-\td g_\alpha(X_\alpha,Z_\alpha)Y_\alpha) +[\td\sigma^\alpha_{\alpha\alpha}(X_\alpha), \td\sigma^\alpha_{\alpha\alpha}(Y_\alpha)]Z_\alpha,
$$
where $X_\alpha,Y_\alpha,Z_\alpha\in {\mathfrak m}_\alpha$ and $R_{\td g_\alpha}$ is the curvature tensor of the metric $\td g_\alpha$.

On the other hand, up to an orthogonal transformation on ${\mathcal H}\subset {\mathfrak m}_0\equiv\bbr^{n_0}$, vectors $H_1,\cdots,H_s$ are uniquely given by the matrix equality
$$
\lmx g(H_1,H_1)&g(H_1,H_2)&\cdots&g(H_1,H_s)\\
g(H_2,H_1)&g(H_2,H_2)&\cdots&g(H_2,H_s)\\
\cdots&\cdots&\cdots&\cdots\\
g(H_s,H_1)&g(H_s,H_2)&\cdots&g(H_s,H_s)\rmx
=
\lmx \td c^2_1&L_1&\cdots&L_1\\
L_1&\td c^2_2&\cdots&L_1\\
\cdots&\cdots&\cdots&\cdots\\
L_1&L_1&\cdots&\td c^2_s\rmx.
$$

Furthermore, it is easily seen from \eqref{000}, \eqref{alp alp 0}, \eqref{alp 0 alp} and \eqref{barsigma0h} that
$\td\sigma^0_{\alpha\alpha}$, $\td\sigma^\alpha_{0\alpha}$, $\td\sigma^\alpha_{\alpha0}$, $\bar\sigma^{\mathcal H}_0$, $\left.\td\sigma_0^{\mathcal H^\bot}\right|_{{\mathcal H^\bot}\times {\mathcal H}}$, $\left.\td\sigma_0^{\mathcal H^\bot}\right|_{{\mathcal H}\times {\mathcal H^\bot}}$ and $\left.\td\sigma_0^{\mathcal H}\right|_{{\mathcal H}\times {\mathcal H}}$ are complete determined by the flat metric $g$, the vectors $H_1,\cdots,H_s$ and the affine mean curvature $L_1$.

Finally, since $\bar\sigma^{\mathcal H^\bot}_0\in {\mathcal S}_{\bbr^{r-1}}\left(-\fr{(n+1)L_1}r\right)$, it can be realized as the second fundamental form of a flat, minimal Lagrangian submanifolds in the complex projective space $\bbc P^{r-1}\left(-4\fr{(n+1)L_1}r\right)$ of holomorphic curvature $-\fr{4(n+1)L_1}r$. Now a theorem of A-M. Li and G. S. Zhao in \cite{li-zhao94} assures that any of such flat, minimal Lagrangian submanifolds is unique up to holomorphic isometries. Thus $\bar\sigma^{\mathcal H^\bot}_0$ is also unique up to isometries on $\bbr^{r-1}$. It then follows that $\td\sigma^0_{00}$ is also completely determined by the flat metric $g$, the vectors $H_1,\cdots,H_s$ and the affine mean curvature $L_1$ up to isometries on $\bbr^{n_0}$.

Summing up, we have proved the conclusion of Lemma \ref{lem4.3}.\endproof

Now we return to the proof of Theorem \ref{chara}.

Let $C$ be give by the second formula in \eqref{newl1c}. Suitably choosing the constants $c_a (1\leq a\leq r)$, $c_{r+\alpha},\, \lalp (1\leq\alpha\leq s)$, we can also assume the first equality. For each $\alpha=1,\cdots,s$, fix one Riemannian metric $$\td g^{(\alpha)}=\fr{(n+1)L_1}{(n_\alpha+1)\lalp}\td g_\alpha$$ on $\td M_\alpha$. Then by \eqref{alp alp alp}, $$\td\sigma^\alpha_{\alpha\alpha}\in {\mathcal S}_{(\td M_\alpha,\td g^{(\alpha)})}\left(\fr{(n_\alpha+1)(-L_1+\td c_\alpha^2)} {(n+1)L_1}\lalp\right).$$

We claim that
\be\label{claim}
\fr{(n_\alpha+1)(-L_1+\td c_\alpha^2)}{(n+1)L_1}=-1,\mb{\ or equivalently,\ }\td c^2_\alpha=\fr{n-n_\alpha}{n_\alpha+1}(-L_1).
\ee

In fact, multiplying $h_{\alpha\gamma}$ to the both sides of \eqref{*} and then taking sum over $\alpha$ we have
\be\label{**}
1=\sum_{\alpha,\beta}h^{\alpha\beta}h_{\alpha\gamma}=-\fr{n_\alpha+1}{rL_1}h_{\alpha\gamma}.
\ee
Since, by \eqref{alp alp 0}, $h_{\gamma\gamma}=\td c^2_\gamma$ and $h_{\alpha\gamma}=L_1$ for $\alpha\neq\gamma$, the right hand side of \eqref{**}
\begin{align}
-\fr{n_\alpha+1}{rL_1}h_{\alpha\gamma}=&-\fr{n_\gamma+1}{rL_1}\td c^2_\gamma-\sum_{\alpha\neq\gamma}\fr{n_\alpha+1}{rL_1}L_1 =-\fr{n_\gamma+1}{rL_1}\td c^2_\gamma-\fr1r\sum_{\alpha\neq\gamma}(n_\alpha+1)\nnm\\
 =&-\fr{n_\gamma+1}{rL_1}\td c^2_\gamma+\fr1r(n_\gamma+1)-\fr1r\sum_\alpha (n_\alpha+1)\nnm\\
 =&-\fr{n_\gamma+1}{rL_1}\td c^2_\gamma+\fr1r(n_\gamma+1)-\fr1r(n-r+1)\nnm\\
 =&-\fr{n_\gamma+1}{rL_1}\td c^2_\gamma-\fr1r(n-n_\gamma)+1\label{***}
\end{align}
From \eqref{**} and \eqref{***} we easily prove the claim \eqref{claim}.

Now the equality \eqref{claim} shows that $\td\sigma^\alpha_{\alpha\alpha}\in {\mathcal S}_{(\td M_\alpha,\td g^{(\alpha)})}(-\!\!\!\!\lalp)$. It follows from Lemma \ref{nai81} that there exists a parallel and minimal Lagrangian submanifold
$\td x_\alpha:(\td M_\alpha,\td g^{(\alpha)})\to \bbc P^{n_\alpha}(-4\!\!\!\!\lalp)$,
which corresponds to a hyperbolic affine hypersphere $x_\alpha:M^{n_\alpha}_\alpha\to \bbr^{n_\alpha+1}$ with the Berwald-Blaschke metric $\stx{\alpha}{g}$, the affine mean curvature $\lalp$ and the parallel Fubini-Pick form $\stx{\alpha}{A}$ where $(M^{n_\alpha}_\alpha,\stx{\alpha}{g})$ is the noncompact symmetric space dual to $(\td M_\alpha,\td g^{(\alpha)})$. Now consider the Calabi composition $\bar x$ of $r$ points and the $s$ hyperbolic affine hyperspheres $x_\alpha$, with the constants $c_a,c_{r+\alpha}$ chosen previously. Suitably choose the parameters $t^1,\cdots,t^{K-1}$, $K=r+s$, one can arrive at $\td g_0=\sum_\lambda\fr{f_{\lambda+1}C}{(n_{\lambda+1}+1)f_\lambda}(d t^\lambda)^2$. Then by Corollary \ref{corr0} and Lemma \ref{lem4.3} we easily find that the parallel and minimal Lagrangian submanifold corresponding to the hyperbolic affine hypersphere $x$ is isometrically equivalent to $\td x$ since they have the same metric and second fundamental form. It then follows that the original hyperbolic affine hypersphere $x$ is equiaffine equivalent to the above Calabi composition $\bar x$.
\endproof

\section{Classification of locally strongly convex hypersurfaces with parallel Fubini-Pick form --- revisted}

In this section, we use the previous correspondence theorem (Theorem \ref{corr thm} and Theorem \ref{chara} to give an alternative proof for the classification of the locally and strongly convex hypersurfaces with parallel Fubini-Pick form. We should remark that this classification has been proved by Z.J. Hu, H.Z. Li and L. Francken in a totally different way (see \cite{hu-li-vra11}).

First we state the classification theorem as follows:
{\thm\label{cla thm}\text{(cf. \cite{hu-li-vra11})}
Let $x:M^n\to \bbr^{n+1}$ ($n\geq 2$) be a locally strongly convex affine hypersurface with parallel Fubini-Pick form $A$. Then either of the following two cases holds:

$(1)$ With the Berwald-Blaschke metric $g$, the Riemannian manifold $(M^n,g)$ is irreducible and $x$ is locally equiaffine equivalent to

$(a)$ one of the three kinds of quadratic affine spheres: Ellipsoid, elliptic paraboloid and hyperboloid; or

$(b)$ the standard embedding of the Riemannian symmetric space ${\rm SL}(m,\bbr)/{\rm SO}(m)$ into $\bbr^{n+1}$ with $n=\fr12m(m+1)-1$, $m\geq 3$;
or

$(c)$ the standard embedding of the Riemannian symmetric space ${\rm SL}(m,\bbc)/{\rm SU}(m)$ into $\bbr^{n+1}$ with $n=m^2-1$, $m\geq 3$; or

$(d)$ the standard embedding of the Riemannian symmetric space ${\rm SU}^*(2m)/{\rm Sp}(m)$ into $\bbr^{n+1}$ with $n=2m^2-m-1$, $m\geq 3$; or

$(e)$ the standard embedding of the Riemannian symmetric space ${\rm E}_{6(-26)}/{\rm F}_4$ into $\bbr^{27}$.

$(2)$ $(M^n,g)$ is reducible and $x$ is locally affine equivalent to the Calabi product of $r$ points and $s$ of the above irreducible hyperbolic affine spheres of lower dimensions, where $r$, $s$ are nonnegative integers and $r+s\geq 2$.}

\proof First note that when $A\equiv 0$, the Pick invariant $J$ vanishes identically, $x$ must be locally equiaffine equivalent to one of the quadratic affine sphere as mentioned above (\cite{li-sim-zhao93}). Therefore it suffices to consider the case that $A\neq 0$. But if it is the case, then by \cite{dil-vra-yap94}, $x$ must be a hyperbolic affine sphere and thus its affine mean curvature $L_1<0$.

Now let $x:M^n\to\bbr^{n+1}$ be a locally strongly convex hypersurface with parallel Fubini-Pick form $A\neq 0$. Denote by $g$, $L_1$ the Berwald-Blaschke metric on $M^n$ and the affine mean curvature of $x$, respectively. Then by Proposition \ref{sym para1}, the immersion $x$ is locally symmetric and thus $(M^n,g)$ is locally isometric to a simply connected symmetric space $G/K$, still denoted by $M^n$. Let $(\td M^n,\td g)\equiv \td G/K$ be the Riemannian symmetric space dual to $(M^n,g)$. Then $(\td M^n,\td g)$ is also simply connected and, by Theorem \ref{corr thm}, $x$ uniquely defines a symmetric and minimal Lagrangian immersion $\td x:\td M^n\to \bbc P^n(-4L_1)$ of $(\td M^n,\td g)$ into the $n$-dimensional complex space $\bbc P^n(-4L_1)$ with constant holomorphic sectional curvature $-4L_1$.

Now Proposition \ref{sym para2} implies that $\td x$ is also parallel, i.e., the second fundamental form $\sigma$ is parallel.
The symmetric trilinear form $\td\sigma$ given by the second fundamental form $\sigma$ of $\td x$ naturally corresponds to $A$. In particular, $\td\sigma\neq 0$. If $x$ is not a Calabi composition of some points and some hyperbolic affine hyperspheres, then by Theorem \ref{chara}, $(M^n,g)$ is irreducible as a Riemannian manifold. It follows that the Remannian manifold $(\td M^n,\td g)$ is also irreducible. Now we can use one of the main results of Naitoh in \cite{nai81} (Theorem 4.5) to conclude that $(\td M^n,\td g)$ can only be one of the following four:
\be\label{symspaces}\left\{\aligned
&{\rm SU}(m)/{\rm SO}(m),\, n=\fr12m(m+1)-1,\quad m\geq 3;\\
&{\rm SU}(2m)/{\rm SP}(m),\, n=2m^2-m-1,\quad m\geq 3;\\
&{\rm SU}(m)\equiv ({\rm SU}(m)\times {\rm SU}(m))/{\rm SU}(m), n=m^2-1,\quad m\geq 3;
\\
&{\rm E}_6/{\rm F}_4,\,n=26,\endaligned
\right.
\ee
where ${\rm E}_6$ is the compact real Lie group with the unique compact real Lie algebra of type $\mathfrak{e}_6$. In fact, ${\rm E}_6$ can be defined as the Lie group of all the complex linear automorphisms on the complex Jordan algebra $\mathfrak{J}^{\bbc}$, the complexification of the real Jordan algebra
$\mathfrak{J}$, keeping invariant both the determinant function and the standard Hermitian inner product on $\mathfrak{J}^{\bbc}$, while ${\rm F}_4$ is the subgroup of ${\rm E}_6$ of those isomorphisms also keeping the standard inner product invariant, or equivalently, ${\rm F}_4={\rm E}_6\cap {\rm O}(\mathfrak{J}^{\bbc},\bbr)$, and can be identified with the Lie group of all the real Jordan algebra automorphisms on $\mathfrak{J}$.

For each of the above symmetric spaces, Naitoh defines in \cite{nai81} one standard minimal and parallel Lagrangian imbedding into the complex space $\bbc P^n(c)$. Those imbeddings are equivariant and uniquely determined (as parallel Lagrangian immersions) by the constant $c$ up to holomorphic isometries (See \cite{nai81}: Theorems 3.4 and 3.6; Lemma 4.2; Proposition 4.4).

Hence, by the theory of Riemannian symmetric spaces, we know that the Riemannian manifold $(M^n,g)$ must be one of the following four spaces dual to those in \eqref{symspaces}:
\be\label{dualspaces}\left\{\aligned
&{\rm SL}(m,\bbr)/{\rm SO}(m),\, n=\fr12m(m+1)-1,\quad m\geq 3;\\
&{\rm SU}^*(2m)/{\rm SP}(m),\, n=2m^2-m-1,\quad m\geq 3;\\
&{\rm SL}(m,\bbc)/{\rm SU}(m), n=m^2-1,\quad m\geq 3;\\
&{\rm E}_{6(-26)}/{\rm F}_4,\,n=26,
\endaligned\right.
\ee
where ${\rm SU}^*(2m)={\rm SL}(2m,\bbc)\cap {\rm U}^*(2m)$ with ${\rm U}^*(2m)$ the usual ${\rm U}$-star group of order $2m$, and ${\rm E}_{6(-26)}$ is one of the noncompact real forms of type $\mathfrak{e}_6$ with ${\rm F}_4$ as its maximal compact subgroup. As a matter of fact, ${\rm E}_{6(-26)}$ is defined as the Lie group of all real linear automorphisms on the $27$-dimensional real Jordan algebra
$\mathfrak{J}$ that keeps the standard determinant function invariant, and ${\rm F}_4$ is identified with the Lie group of all the elements in ${\rm E}_{6(-26)}$ that keeps invariant the identity matrix $I_3\in\mathfrak{J}$.

We would like to interrupt here to remark that all of these four symmetric spaces has been discussed in \cite{sas80} where the author proved these four examples are homogeneous hyperbolic affine spheres. On the other hand, for each of the first three symmetric spaces listed above, there can be defined one equivariant imbedding of it into $\bbr^{n+1}$ with a given affine mean curvature $L_1$. They are naturally symmetric affine hypersurfaces and, by Proposition \ref{sym para1}, are of parallel Fubini-Pick forms. In fact, the standard equivariant imbedding of ${\rm SL}(m,\bbr)/{\rm SO}(m)$ into $\bbr^{\fr12m(m+1)}$ can be found in \cite{nom-sas94}; The same idea was used by O. Birembaux and M. Djoric in \cite{bir-djo12} to define and study the equivariant standard imbeddings of ${\rm SU}^*(2m)/{\rm SP}(m)$ and ${\rm SL}(m,\bbc)/{\rm SU}(m)$ with the special case $m=3$. For general $m$, these examples can be found in \cite{hu-li-vra11}. We also remark that in the same paper \cite{bir-djo12}, the authors also defined an explicit imbedding of ${\rm E}_{6(-26)}/{\rm F}_4$ into $\bbr^{27}$ in terms of coodinates which is proved to have parallel Fubini-Pick form by Hu etc in \cite{hu-li-vra11}. But here we would like to show that the same idea used in \cite{nom-sas94} and rather recently by \cite{bir-djo12} can be extended to define an affine equivariant imbedding of $M:={\rm E}_{6(-26)}/{\rm F}_4$ into $\bbr^{27}$. This treatment is more natural one and has not appeared in the literatures published. The details are presented as follows:

Let $\mathbb{O}$ be the space of octonions and $\mathfrak{J}$ be the set of $3\times 3$ Hermitian matrices with entries in $\mathbb{O}$, that is
$$
\mathfrak{J}=\{X=\lmx \xi_1&x_3&\bar x_2\\ \bar x_3&\xi_2&x_1\\
x_2&\bar x_1&\xi_3\rmx\in {\rm M}(3,\mathbb{O});\ \bar X^t=X\},
$$
where ${\rm M}(3,\mathbb{O})$ is the real vector space of all octonian square matrices of order $3$. Clearly $\mathfrak{J}$ is a real vector space of dimension $n+1:=27$ and thus can be identified with $\bbr^{27}$. On $\mathfrak{J}$, the symmetric Jordan multiplication $\circ$ and the standard inner product $(\cdot,\cdot)$ on $\mathfrak{J}$ are defined as follows:
$$
X\circ Y=\fr12(XY+YX),\quad (X, Y)=\tr(X\circ Y).
$$
Furthermore, the cross product $\times$ and the determinant function $\det$ are given by
\bea
&X\times Y=\fr12(2X\circ Y-\tr(X)Y-\tr(Y)X+(\tr(X)\tr(Y)-\tr(X\circ Y))I_3)\\
&\det(X)=\fr13(X\times X, X).
\eea

The noncompact group ${\rm E}_{6(-26)}$ is defined as the set of all determinant-preserving real linear automorphism on $\mathfrak{J}$, that is
\be\label{E6(-26)} {\rm E}_{6(-26)}=\{A\in {\rm GL}_\bbr(\mathfrak{J});\ \det(AX)=\det(X),\,\forall X\in\mathfrak{J}\}.
\ee
The maximal compact subgroup of ${\rm E}_{6(-26)}$ is given by
\bea
&&\hs{-1.1cm}{\rm F}_4=\{A\in {\rm E}_{6(-26)};\ A(X\circ Y)=(AX)\circ(AY),\, \forall X,Y\in\mathfrak{J}\}\label{F4-1}\\
&&\hs{-.6cm}\equiv\{A\in {\rm E}_{6(-26)};\ A(I_3)=I_3\}.\label{F4-2}
\eea
\vs{.3cm}
For each matrix $T\in\mathfrak{J}$, there associated an element $\td T\in {\rm E}_{6(-26)}$ defined by
$$
\td T(X):=T\circ X,\quad \forall X\in \mathfrak{J}.
$$
Define
$$
\mathfrak{m}=\{\td T;\ T\in\mathfrak{J}_0\}, \text{\ where\ }
\mathfrak{J}_0=\{T\in\mathfrak{J};\ \tr T=0\}.
$$
Denote by ${\mathfrak f}_4$ the Lie algebra of ${\rm F}_4$. Then by \cite{yok09}, the Lie algebra $\mathfrak{e}_{6(-26)}$ has a canonical vector space decomposition as
\be\label{dec1}\mathfrak{e}_{6(-26)}=\mathfrak{f}_4+\mathfrak{m}\ee
satisfying
$[\mathfrak{f}_4,\mathfrak{m}]\subset \mathfrak{m}$, $[\mathfrak{m},\mathfrak{m}]\subset \mathfrak{f}_4$. Note that we have a natural identification $\mathfrak{m}\equiv T_oM$ where $o:={}_{I_{27}}F_4$ with $I_{27}$ the identity element in $E_{6(-26)}$.

{\prop\label{prop 5.1}
${\rm E}_{6(-26)}$ is a subgroup of the special linear group ${\rm SL}(27,\bbr)$.}

\proof To prove Proposition \ref{prop 5.1}, we first define
$$
E_1=\lmx 1&0&0\\0&0&0\\0&0&0\rmx, \, E_2=\lmx 0&0&0\\0&1&0\\0&0&0\rmx,\,
E_3=\lmx 0&0&0\\0&0&0\\0&0&1\rmx;
$$
$$
F_1(x)=\lmx 0&0&0\\0&0&x\\0&\bar x&0\rmx,\,F_2(x)=\lmx 0&0&x\\0&0&0\\ \bar x&0&0\rmx,
\,F_3(x)=\lmx 0&x&0\\ \bar x&0&0\\0&0&0\rmx,\,x\in\mathbb{O}.
$$
Then $\{E_i,F_i(x); x\in\mathbb{O}, i=1,2,3\}$ generates $\mathfrak{J}$ and for $x,y\in\mathbb{O}$
\bea
&&\hs{-.6cm}\left\{\aligned &E_i\circ E_i=E_i,\\&E_i\circ F_i(x)=0,\\&F_i(x)\circ F_i(y)=(x,y)(E_{i+1}+E_{i+2}),\endaligned\right.\label{ef1}\\
&&\hs{-.6cm}\left\{\aligned &E_i\circ E_j=0,\,i\neq j,\\&E_i\circ F_j(x)=\fr12 F_j(x),\,i\neq j,\\
&F_i(x)\circ F_{i+1}(y)=\fr12F_{i+2}(\ol{xy}),\endaligned\right.\label{ef2}
\eea
where the indices are considered as being modulo to $3$.

Furthermore, we define
\bea
&\mathfrak{M}^-=\{A\in {\rm M}(3,\mathbb{O});\ \bar A^t+A=0\},\\
&\mathfrak{M}^-_0=\{A\in \mathfrak{M}^-;\ {\rm diag}\,A=0\};\\
&\mathfrak{der}=\{\delta\in\mathfrak{f}_4;\ \delta E_1=\delta E_2=\delta E_3=0\}\eea
where ${\rm diag}\,A=0$ means that all the diagonal elements of $A$ vanish. On ${\rm M}(3,\mathbb{O})$, there is a natural bracket $[\cdot,\cdot]$ given by
$$
[X,Y]:=XY-YX.
$$
Then we have

{\lem\label{lem 4.2} {\rm(\cite{yok09}, Lemma 2.3.3)}
$$[\mathfrak{M}^-,\mathfrak{J}]\subset \mathfrak{J},\ [\mathfrak{J},\mathfrak{J}]\subset \mathfrak{M}^-.
$$}
Clearly, for each $A\in \mathfrak{M}^-$, there associated one element $\td A\in \mathfrak{f}_4$ such that
\be\label{tdA}\td A(X):=[A,X],\quad \forall X\in\mathfrak{J}.\ee
On the other hand, for each $\delta\in\mathfrak{der}$, there corresponds uniquely one $D_1\in \mathfrak{so}(\mathbb{O})$ such that (\cite{yok09}, Proposition 2.3.7)
\be\label{delta}
\delta\lmx \xi_1&x_3&\bar x_2\\ \bar x_3&\xi_2&x_1\\ x_2&\bar x_1&\xi_3\rmx
=\lmx 0&D_3x_3&\ol{D_2x_2}\\ \ol{D_3x_3}&0&D_1x_1\\ D_2x_2&\ol{D_1x_1}&D_3\xi_3\rmx
\ee
where $D_2,D_3\in\mathfrak{so}(\mathbb{O})$ are determined by the following equation
$$
(D_1x)y+x(D_2y)=\ol{D_3(\ol{xy})},\quad\forall\,x,y\in\mathbb{O}.
$$
If we denote $\td{\mathfrak{M}}^-_0=\{\td A;\ A\in\mathfrak{M}^-_0\}$, then
the Lie algebra
$\mathfrak{f}_4$ of ${\rm F}_4$ is decomposed further into
\be\label{dec2}
\mathfrak{f}_4=\mathfrak{der}+\td{\mathfrak{M}}^-_0.
\ee
Note that $\mathfrak{J}$ is generated by $\{E_i,F_i(x);\ x\in\mathbb{O}\}$ from which one can find an orthonormal basis for  $\mathfrak{J}$. A direct computation in terms of this basis by using \eqref{ef1}, \eqref{ef2}, \eqref{tdA}, \eqref{delta} and the decompositions \eqref{dec1}, \eqref{dec2} shows that

{\lem\label{lem 4.3}
For each $A\in \mathfrak{e}_{6(-26)}$, $\tr A=0$.}

Then Proposition \ref{prop 5.1} follows immediately. \endproof

{\expl The equivariant imbedding of ${\rm E}_{6(-26)}/{\rm F}_4$}

For any given constant $L_1<0$, set
$$C=\sqrt{3}(-3L_1)^{-\fr{n+2}2}>0.$$
Then define a smooth map $f:{\rm E}_{6(-26)}\to \mathfrak{J}$ by
$f(L)=C\cdot L(I_3)$ for all $L\in {\rm E}_{6(-26)}$. Clearly, for any $L_1,L_2\in {\rm E}_{6(-26)}$, $f(L_1)=f(L_2)$ if and only if $(L^{-1}\circ L_2)(I_3)=I_3$. By the definition of ${\rm F}_4$, $f$ naturally induces a smooth map $x:{\rm E}_{6(-26)}/{\rm F}_4\to\bbr^{27}\equiv\mathfrak{J}$:
\be\label{e6/f4}
x({}_L{\rm F}_4)=C\cdot L(I_3),\quad \forall L\in {\rm E}_{6(-26)}.
\ee

By Proposition \ref{prop 5.1}, we can choose a volume element on $\bbr^{27}$,
say, the canonical volume element with respect to the inner product $(\cdot,\cdot)$ on $\mathfrak{J}$, so that ${\rm E_{6(-26)}}$ can be identified with a subgroup of the
group ${\rm UA}(27)$ of unimodular affine transformation on $\bbr^{n+1}$.
Therefore, the induced map $x$ is equivariant as an affine hypersurface
in $\bbr^{n+1}$. Consequently all the equiaffine invariants of $x$ such as the Berwald-Blaschke metric, the Fubini-Pick form and the fundamental form are ${\rm E_{6(-26)}}$-invariant.

Now for each $\td T\in\mathfrak{m}\equiv T_oM$, $T\in\mathfrak{J}_0$, it holds clearly that
$$
x_*(\td T)=\left.\dd{}{t}\right|_{t=0}(\exp t\td T(I_3))=\td T(I_3)=C(T\circ I_3)=C\cdot T.
$$
This shows that $x$ is an immersion at $o$ and thus is an immersion globally since $x$ is equivariant. Clearly, $x$ is monomorphic and so is an imbedding of $M$ into $\bbr^{27}$.

Moreover, since for each $T\in \mathfrak{J}_0$,
$$(T,I_3)=\tr(T\circ I_3)=\tr T=0,$$ $x(o)$ is a transversal vector of $x$ at $o$ and thus is transversal everywhere. Furthermore,
for all $X,Y\in\mathfrak{J}_0$,
\bea
&&\hs{-1.8cm}\td X(x_*(\td Y))=\left.\dd{}{t}\right|_{t=0}x_*\left((L_{\exp t\td X})_*(\td Y)\right)\nnm\\
&&\hs{-.1cm}=\left.\ppp{}{t}{s}\right|_{t=s=0}(C\exp t\td X\cdot\exp s\td Y(I_3))\nnm\\
&&\hs{-.1cm}=C(X\circ(Y\circ I_3))=C(X\circ Y)\nnm\\
&&\hs{-.1cm}=C\left(X\circ Y-\fr13\tr(X\circ Y)I_3\right)+\fr13C(X,Y)I_3\label{gaussf}
\eea
implying that $x$ is locally strongly convex since $(X,Y)=\tr (X\circ Y)$ is positive definite.

Note that the inner product $(\cdot,\cdot)$ on $\mathfrak{J}_0$ is $\mathfrak{f}_4$-invariant and that the correspondence $\ \widetilde{\ }:\mathfrak{J}_0\to \mathfrak{m}$ is $\mathfrak{f}_4$-equivariant. It follows by the definition (cf. \eqref{dfn h} and \eqref{dfn g}) that the Berwald-Blaschke metric $g$ of $x$ is the invariant metric on ${\rm E}_{6(-26)}/{\rm F}_4$ induced by
$$
g_o(\td X,\td Y):=\left(\fr1{\sqrt{3}}C\right)^{\fr2{n+2}}(X,Y)=-\fr1{3L_1}(X,Y),\quad \forall\, X,Y\in\mathfrak{J}_0,
$$
or, equivalently $g_o(\td X,\td Y)=-\fr1{3L_1}\tr(X\circ Y)$. Taking the trace of \eqref{gaussf} respect to the Berwald-Blaschke metric $g$, we find that the affine normal $\xi=-L_1\cdot x$ at $o$ and thus at everywhere. It follows that $x$ is a hyperbolic affine hypersphere with the affine mean curvature being the given number $L_1$.

On the other hand, the invariant Fubini-Pick form $A$ of $x$ is induced by the following $\mathfrak{f}_4$-invariant form $A_o$:
$$
A_o(\td X,\td Y,\td Z)=g_o\left(-L_1\left(X\circ Y-\fr13\tr(X\circ Y)I_3\right)\widetilde{},\td Z\right),\,
\forall X,Y,Z\in \mathfrak{J}_0.
$$
In particular, $x$ is a symmetric equiaffine sphere in $\bbr^{27}$. Now we use Proposition \ref{sym para1} to conclude that the Fubini-Pick form $A$ is parallel.

Now we come back to complete the proof of Theorem \ref{cla thm}.

By the conclusions of Naitoh (\cite{nai81}: Lemma 4.2 and Proposition 4.4), we know that, for each of the symmetric spaces $\td G/K$ listed in \eqref{symspaces}, the invariant metric $\td g$ and the totally symmetric invariant form $\td\sigma$ on $\td G/K$ are uniquely determined by the constant $-L_1$, or equivalently via Theorem \ref{corr thm}, for each of the spaces $G/K$ listed in \eqref{dualspaces}, the invariant Berwald-Blaschke metric $g$ and the invariant Fubini-Pick form $A$ are uniquely determined by the affine mean curvature $L_1$. It then follows from Theorem \ref{affine uniqueness} that the given irreducible affine hypersurface $x:M^n\to\bbr^{n+1}$ in Theorem \ref{cla thm} with parallel Fubini-Pick form is affine equivalent to one of the standard imbeddings of the spaces in \eqref{dualspaces} in $\bbr^{n+1}$.\endproof

\flushleft
Xingxiao Li\\
School of Mathematics and Information Sciences\\
Henan Normal University\\
XinXiang Henan 453007\\
P.R.China\\
email: xxl@henannu.edu.cn

\end{document}